\documentclass{siamltex}

\usepackage{amsmath}
\usepackage{amssymb}

\usepackage[toc,page]{appendix}

\usepackage{pgfplots}
\usepgfplotslibrary{groupplots}
\usepackage[ComplainAboutUnexternalized]{externaltikz}
\usepackage{calc}

\usepackage{subcaption}
\usepackage{hyperref}
 \usepackage{mathtools}

\newtheorem{remark}{Remark}

\begin{document}
\title{A kinetic traffic network  model and its macroscopic limit: diverging lanes}

\author{R. Borsche\footnotemark[1] 
	\and  A. Klar\footnotemark[1] \footnotemark[2]}
\footnotetext[1]{Technische Universit\"at Kaiserslautern, Department of Mathematics, Erwin-Schr\"odinger-Stra{\ss}e, 67663 Kaiserslautern, Germany 
	(\{borsche, klar\}@mathematik.uni-kl.de)}
\footnotetext[2]{Fraunhofer ITWM, Fraunhoferplatz 1, 67663 Kaiserslautern, Germany}

	\maketitle

	\begin{abstract}
		This paper is a continuation of the work in \cite{BK20} on macroscopic limits of kinetic traffic networks.
		There, the case of merging lanes has been investigated. In  the present work we propose coupling conditions for a kinetic two velocity model for vehicular traffic for  junctions with diverging lanes.
		We consider cases with and without directional preferences and present corresponding kinetic coupling conditions.
		From this kinetic network model coupling conditions for a macroscopic traffic model are  derived.
		We use as in \cite{BK20}  an  analysis of  the layer equations at the junction in combination with a suitable matching procedure with half-Riemann problems for the macroscopic model. In this way  classical coupling conditions for scalar conservation laws for traffic flow on networks are derived from an underlying network problem.
	\end{abstract}

\section{Introduction}

In the context of traffic flow on networks many models rely on hyperbolic partial differential equations ranging  from scalar conservation laws over systems of conservation laws  to kinetic models \cite{GPBook,AR,AKMR,PSTV17,PSTV171,HPRV20}.
For some of these equations hierarchies of such models have been established for single roads, e.g. deriving macroscopic equations from microscopic ones. Such hierarchies  have been investigated for example  in \cite{AR,AKMR,Ber,Dag2}.
Similarly macroscopic models can be obtained from kinetic descriptions, see  e.g.  \cite{Hel,KW97}.
On the other hand there has been a continuous effort in extending such models onto networks of roads. 
Most of theses approaches consider scalar conservation laws, see e.g. \cite{HR95,CGP05,GPBook,CCDT19}.
Only few attempts have been made for second order or kinetic problems \cite{HM09,HR}.
In none of these works a hierarchy of such models on a network has been investigated and  network models 
 for macroscopic models derived  from underlying kinetic or microscopic models.
 
 In \cite{BK20} such an investigation has been  started for a basic kinetic model leading to scalar hyperbolic traffic models.
 There, we have only considered the case of merging junctions.
The present paper aims to close this gap for diverging junctions, i.e. for junctions with one ingoing and 
two outgoing roads.
Hereby we follow closely the strategy developed in \cite{BKKP16,BK18b,BK18c,BK20}.
In \cite{BK18c} coupling conditions for athe wave equation  have been derived from an underlying linear kinetic description.
For the wave equation this procedure involves an approximation of the layers arising at the junction \cite{BK18c}.
Non linear problems, like the Burgers equation, further require a half Riemann problem to link the layer solution to the macroscopic states \cite{BK18b}.
By a successive combination of these tools coupling conditions for the associated macroscopic problems can be derived from the kinetic models.

In the present paper we first propose coupling conditions on the kinetic level for a two velocity traffic model, derived in \cite{BK18}.
An important advantage, compared to the macroscopic equation, is ,that for the two velocity model the required number of coupling condition remains constant.
On each individual road the equations contain a scaling parameter $\epsilon$ and their solutions converge for $\epsilon \rightarrow 0$ to an associated scalar traffic flow model.
If $\epsilon$ is send to zero on the network, boundary layers at the junctions can arise.
The structure of such layers can be studied by investigating the solutions of the associated half space problems.
Similar approaches have been used in \cite{BSS84,BLP79,G08,CGS,UTY03,N99,AM04} in the context of kinetic equations and in \cite{WY99,WX99,LX96,X04} for hyperbolic relaxation systems.
Since the equations under consideration are non linear, these half space problems have to be coupled to half Riemann problems of the macroscopic equation, as outlined in \cite{BK18b}.
Combining all these, a non linear problem is established in terms of the macroscopic unknowns at the node.
The solution to this coupling problem provides the required states at the junction for the scalar hyperbolic model
and the limit of the kinetic traffic network problem  as  $\epsilon\rightarrow 0$.

The  paper is organized as follows.
Fist the two velocity traffic model is revisited and, shortly, the main features are collected.
In section \ref{sec:kineticcouplingconditions} coupling conditions for this kinetic equation are proposed
for diverging junctions.
The case of a junction with one ingoing and two outgoing roads is considered and  drivers with or without directional preferences are discussed.
Aside the conservation of mass at the node, the free space available in the exiting roads is an important quantity.
In section \ref{macroscopiccc} the resulting coupling conditions for the scalar variables as $\epsilon \rightarrow 0$
are presented. 
A detailed derivation of these results is given in section \ref{technical}.
Further, the obtained results are verified by numerical examples in section \ref{Numerical results}.

\section{Kinetic and macroscopic traffic equations}
\label{equations}

We consider a minimal kinetic discrete velocity model \cite{BK18} with just two velocities $v_0 = 0$ and $v_1=1$.
The densities corresponding to these velocities are $f_0$, which represents the  stopped cars, and $f_1$, which is the density of driving cars. 
Using these we can define the total density of cars as $\rho = f_0+f_1\in [0,1]$ and the mean flux as $q =v_0 f_0 +v_1 f_1= f_1 $ or reversely 
$$
f_0 =\rho-q\ ,\quad  f_1 =q\ .
$$
The dynamics of these quantities is governed by the discrete velocity model developed in \cite{BK18}
\begin{eqnarray}
\label{eq:dvm01}
\begin{aligned}
\partial_t f_0  - \frac{1}{1-\rho }    f_1\partial_x f_0&= -\frac{1}{\epsilon} \left(f_0-\rho +F(\rho) \right)  \\
\partial_t f_1  + \partial_x f_1 + \frac{1}{1-\rho }   f_1\partial_x f_0&=  -\frac{1}{\epsilon} \left(f_1- F(\rho) \right)\ ,
\end{aligned}
\end{eqnarray}
where $F= F(\rho)$ is a given traffic density-flow function or fundamental diagram. 
We assume $F:[0,1]\rightarrow[0,1]$ to be a smooth function with $F(0)=0=F(1)$,  $F^\prime(\rho) \le 1$ and its graph in the triangle $0 \le \rho \le 1, 0 \le q \le \rho$.
In the following we restrict ourselves to strictly concave fundamental diagrams $F$.
The point, where the maximum of $F$ is attained we denote by $\rho^{\star}$ and the maximal value by $ F(\rho^{\star})=\sigma$.

The two eigenvalues corresponding to \eqref{eq:dvm01} are $\lambda_1 = - \frac{q}{1-\rho}\leq 0 <\lambda_2 = 1$ with the respective eigenvectors
$r_1 = \left( 1,\lambda_1 \right)^T, r_2 = \left(1,1\right)$. 
The system is strictly hyperbolic and both characteristic families are linearly degenerate.
The integral curves of the hyperbolic system are given by $q= q_L \frac{1-\rho}{1-\rho_L}$ for the 1-field and by $q = \rho -\rho_R+q_R$ for the 2-field. 
As the maximal velocity is set to $1$ the region $0 \le \rho \le 1, 0 \le q \le \rho$ is an invariant region for the kinetic equations.

A Riemann invariant of the first characteristic family is 
$$Z= \frac{q}{1-\rho+q} = \frac{q}{1-w}= \frac{f_1}{1-f_0}\in [0,1]\ .$$
As for $q$ we have $ 0 \le Z \le \rho$.
By using this variable $Z$ the system \eqref{eq:dvm01} can be transformed into conservative form
	\begin{align}\label{eq:lindeg+relax}
	\begin{aligned}
	\partial_t \rho + \partial_x q&=0\\
	\partial_t Z + \partial_x Z &= -\frac{(1-Z)}{\epsilon (1-\rho)} \left(q - F(\rho)  \right)
	\end{aligned}
	\end{align}
with $q =  \frac{Z}{1-Z}(1-\rho)$.
Note that this change does not influence the speed of possible discontinuities, as both fields are linearly degenerate.

A  Riemann invariant of the  second characteristic family is 
$$w =\rho - q = f_0 \in [0,1]\ .$$ 
the quantity $1-w$ can be interpreted as the free space available or the maximal possible number of driving cars.
$Z$ can be understood as the ratio between the actual number of driving cars and the maximal possible number of driving cars.
	
Similarly, equation \eqref{eq:dvm01} can be expressed in macroscopic variables $\rho$, $q$ as
\begin{align}\label{macro0}
\begin{aligned}
\partial_t \rho + \partial_x q &=0\\
\partial_t q + \frac{ q}{1-\rho}  \partial_x \rho  + (1-\frac{ q}{1-\rho})\partial_x  q  &=-\frac{1}{\epsilon} \left(q-F(\rho) \right) \ .
\end{aligned}
\end{align}
Concerning the convergence of its solutions towards the solutions of the scalar conservation law $\partial_t \rho + \partial_x F(\rho) =0$ as $\epsilon$ tends to $0$ the subcharacteristic condition has to be satisfied \cite{LX96}.
Setting $q = F(\rho)$ in the eigenvalues, the subcharacteristic condition states
$$
-  \frac{ F(\rho)}{1-\rho} \le F^\prime(\rho) \le 1 \  \mbox{ for } \  0 \le \rho \le 1\ .
$$
\begin{remark}
	The condition is fulfilled for strictly concave fundamental diagrams $F$.
	For example, in the classical LWR case with $F(\rho) =\rho (1-\rho) $ and  $F^\prime (\rho) = 1- 2 \rho$ the above  condition is
	$$
	-   \rho \le 1- 2 \rho \le 1\  \mbox{ for } \ 0 \le \rho \le 1\ ,
	$$
	which is obviously satisfied. 
\end{remark}

On a finite domain the kinetic problem \eqref{eq:dvm01} has to equipped with boundary conditions.
At the left boundary at $x=x_L$ a value for the 2- Riemann invariant $Z(x_L) = \frac{q(x_L)}{1-\rho(x_L)-q(x_L)} =  \frac{f_1(x_L)}{1-f_0(x_L)}$ is set and for the right boundary $x= x_R$ the 1-Riemann invariant $w(x_R)=f_0(x_R)$.
As the first eigenvalue is always non-positive and the second is constant, the number of boundary conditions does not alter.

\section{Kinetic Coupling conditions}
\label{sec:kineticcouplingconditions}
In this section we propose coupling conditions for the kinetic two-velocity model \eqref{eq:dvm01}.
As on each road there is exactly one outgoing characteristic family, we have to provide  three conditions  at a junction connecting three roads. 
In any case the conservation of mass will be imposed, i.e. all cars entering a junction via one of the incoming roads will exit on one of the outgoing roads.
For the remaining two conditions we will analyze the maximal possible number of driving cars $1-w = 1-f_0$ on the exiting roads.  
Note that only the stopped cars $w = f_0$  on the outgoing roads can block the traffic, as the driving ones $f_1=q$ will clear the space.
Thus $1-f_0$ is the available free space for driving cars.

From the mathematical point of view it is important to supply values for the correct characteristic variables.
We denote with $\hat{\cdot}$ the known traces at the junction. 
The unknown characteristic variables as well as partially known quantities in the junctions do not have a superscript.

\begin{figure}[h!]
	\begin{center}
		\externaltikz{sketch_12node}{
			\begin{tikzpicture}[thick]
			\def\len{2}
			\node[fill,circle] (N) at (0,0){};
			\draw[->] (-\len,0)--(N) node[above,pos = 0.5]{$1$};
			\draw[->] (0,0)--(\len,0.6) node[above,pos = 0.5]{$2$};
			\draw[->] (0,0)--(\len,-0.6) node[below,pos = 0.5]{$3$};
			\end{tikzpicture}
		}
	\end{center}
	\caption{A junction with one incoming and two outgoing roads (1-2 node).}
	\label{fig:split}
\end{figure}
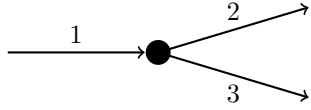

We consider junctions with one incoming and two outgoing roads.
Road number $1$ is the incoming one, see Figure \ref{fig:split}.

\subsection{Junction with diverging lanes and no driver preferences}	
\label{nopref}

In this section we assume that the drivers have no priority which road to take. 
Drivers decide locally according to the free space which road they want to follow. 
Thus the free space in the exiting roads is the relevant quantity. 

The free space observed by the drivers in road $1$ when looking towards road $2$ is the free space in road $2$ plus the space occupied by the cars driving towards  road $3$
\begin{align}\label{eq:split_kin_eq1}
1-w^1= 1-f^1_0 = 1- f_0^2 + f_1^3= 1-w^2+ q^3\ .
\end{align}
Similarly, we obtain  
\begin{align}\label{eq:split_kin_eq2}
1-w^1= 1-f^1_0 = 1- f_0^3 + f_1^2 = 1-w^3+q^2\ .
\end{align}
Finally the conservation of mass is
\begin{align}\label{eq:split_kin_eq3}
f_1^1  = q^1= q^2+q^3=f_1^2 + f_1^3\ .
\end{align}
Note that \eqref{eq:split_kin_eq1} and \eqref{eq:split_kin_eq2} not only specify the free space on road one, but also serve as condition for the distribution of the cars.

It is easy to see that these equations are only valid under further restrictions, e.g. if on road $2$ there are no stopped cars, $f_0^2=w^2=0$, and $f_1^3>0$ equation \eqref{eq:split_kin_eq1} will result in a free space larger than one. 

Thus, these conditions have to be truncated such that all quantities remain within their physical bounds.  
They can be more easily determined expressing the equations   (\ref{eq:split_kin_eq1}-\ref{eq:split_kin_eq3}) in characteristic variables or Riemann Invariants. 
We have
\begin{align*}
1-{w}^1 &= 1-\hat{w}^2+{Z}^3(1-\hat{w}^3)\\
1-{w}^1 &= 1-\hat{w}^3+{Z}^2(1-\hat{w}^2)\\
\hat{Z}^1\left(1-{w}^1\right) &= {Z}^2\left(1-\hat{w}^2\right)+{Z}^3(1-\hat{w}^3).
\end{align*}
The above equations  are solved by
\begin{align*}
	1-w^1&=\frac{2-(\hat w^2+\hat{w}^3)}{2-\hat Z^1}
\end{align*}	
and 
\begin{align*}
(1-\hat w^2)Z^2 = \frac{1}{2-\hat Z^1}(\hat Z^1 (1-\hat w^3)-(\hat w^2-\hat w^3))\\
(1-\hat w^3)Z^3 = \frac{1}{2-\hat Z^1}(\hat Z^1 (1-\hat w^2)+(\hat w^2-\hat w^3))\ .
\end{align*}

We obtain a valid expression for $0 \le 1-w^1 \le 1$ if 	
$
2-(\hat w^2+\hat{w}^3) \le 2-\hat Z^1
$
or
$\hat w^2+\hat{w}^3 \ge \hat Z^1$.
Moreover, if additionally $\hat Z^1 (1-\hat w^3) \ge \hat w^2 - \hat w^3$ and $\hat Z^1 (1-\hat w^2) \ge \hat w^3 - \hat w^2$ we obtain  admissible equations for $0 \le Z^2, Z^3\le 1$.
We note that the last conditions give $q^1 \ge \   \vert \hat w^2- \hat w^3 \vert $ due to the fact that $1-w^1$ is larger than $1-\hat w^2, 1-\hat w^3$.

In this case we have 
\begin{align*}
	(1-\hat{w}^2)+q^3
	=(1-\hat{w}^3)+q^2
	\end{align*}
and therefore using $q^1= q^2+q^3$ we obtain
$$
\rho^2= \rho^3\ .
$$
Moreover,
\begin{align*}
	(1-\rho^1+q^1)(2-\frac{q^1}{1-\rho^1+q^1})&=2-(\rho^2+\rho^3)+q^1\\
\end{align*}
gives
$
\rho^1 = \frac{\rho^2+\rho^3}{2}
$ and therefore the following conditions
\begin{align*}
	\rho^1 = \rho^2=\rho^3\ , 
	q^1=q^2+q^3\ .
\end{align*}		
Using  $1-w^1 = \frac{1}{2}(1-\hat w^2+1-\hat w^3 +q^1)$ a direct computation gives 
\begin{align}
\label{fluxdistribution}
		\begin{aligned}
			q^2 = \frac{q^1}{2} -\frac{\hat w^2-\hat w^3}{2}\\
			q^3 = \frac{q^1}{2} +\frac{\hat w^2-\hat w^3}{2}\ .
		\end{aligned}
	\end{align}		
That means we have a linear distribution of the outgoing fluxes according to the difference in free space on the outgoing roads.

Finally, we have to discuss situations where the above conditions are not fulfilled.

First we consider the case where still, but $\hat w^2+\hat{w}^3 \le \hat Z^1$. 
We note that in this case $q^1 \ge \   \vert \hat w^2- \hat w^3 \vert $ is automatically satisfied, see Figure \ref{fig:nopref}.
In this case, we consider the interface $\hat w^2+\hat{w}^3 = \hat Z^1$, where,  $1-w^1=1$.
Thus, $w^1 =0$ and then $\rho^1 = q^1$, which we choose as the first condition in this domain.

Moreover, we assume in this case that the above distribution of fluxes \eqref{fluxdistribution} is still valid which leads to
\begin{align*}
(1-\hat w^2)Z^2 = \frac{1}{2}(\hat Z^1 -(\hat w^2-\hat w^3))\\
(1-\hat w^3)Z^3 = \frac{1}{2}(\hat Z^1 +(\hat w^2-\hat w^3))\ .
\end{align*}

We note that the balance of fluxes is still guaranteed for this definition.
We obtain directly 
$\rho^2= \rho^3$.

Next, we consider  the case, where $w^2 \ge w^3 $ and  $\hat Z^1 (1-\hat w^3) \le w^2-w^3 $.
Considering again the interface to the first case, we have in this case 
\begin{align*}
1-w^1&=1-\hat w^3
\end{align*}
and the original equations lead to 
\begin{align*}
Z^2 &=0\\
Z^3 &= \frac{(1-w^1)\hat Z^1}{1-\hat w^3}=\hat Z^1\ .
\end{align*}
This yields 
\begin{align*}
q^2=0\ ,\quad 
q^3 =  q^1\ .
\end{align*}
Finally, we consider the case $w^2 \le w^3 $ and  $\hat Z^1 (1-\hat w^2) \le w^3-w^2 $.
This leads to 
\begin{align*}
1-w^1&=1-\hat w^2\\
Z^2 &= \frac{(1-w^1)\hat Z^1}{1-\hat w^2}= \hat Z^1\\
Z^3 &=0
\end{align*}
and equivalently
\begin{align*}
q^2=0\ ,\quad 
q^3 =  q^1\ .
\end{align*}

\begin{figure}[h!]
	\begin{center}
		\externaltikz{Statespace_H1}{
	 	\begin{tikzpicture}[scale = 4]
	 	\def\dr{0.2}
	 	\node[below] at (0,0) {$(0,0)$};
	 	\node[below] at (1,0) {$(1,0)$};
	 	\node[right] at (1,1) {$(1,1)$};
	 	\node[left] at (0,1) {$(0,1)$};
	 	\node[below] at (0.2,0.8) {$\bf IV$};
	 		 	\node[below] at (0.4,0.5) {$\bf I$};
	 		 	\node[right] at (0.03,0.12) {$\bf II$};
	 		 	\node[left] at (0.8,0.2) {$\bf III$};
	 	\node[below] at (0.4,0) {$\hat Z^1$};
	 	\node[left] at (0,0.4) {$\hat Z^1$};
	 	\draw[line width=1pt](0.4,0)--(0,0.4);
	 	\draw[->](0,0)--(1.2,0) node[below]{$\hat w^2$};
	 	\draw[->](0,0)--(0,1.2) node[left]{$\hat w^3$};
	 	\draw(1,0)--(1,1);
	 	\draw(0,1)--(1,1);
	 	\draw[line width=1pt](0.4,0)--(1,1);
	 	\draw[line width=1pt](0,0.4)--(1,1);
	 	\end{tikzpicture}
	 	 }
	\end{center}
	\caption{Domain of definition for coupling conditions without drivers preference.}
	\label{fig:nopref}
\end{figure}

\subsubsection{Summary}	
\label{kincondsummary}
In summary  we have the following  coupling conditions in characteristic variables distinguishing   4 cases, see Figure \ref{fig:nopref}.

{\bf Case I:} $\hat Z^1 (1-\hat w^3)  \ge \hat w^2- \hat w^3 $    and  $\hat Z^1 (1-\hat w^2) \ge  \hat w^3- \hat w^2$ and $\hat w^2 + \hat w ^3  \ge \hat Z^1$:
$$
\rho^1 = \rho^2 =  \rho^3,\  q^1=q^2+q^3\ ,
$$

{\bf Case II:} $\hat w^2 + \hat w ^3  \le \hat Z^1$
$$\rho^1 = q^1,\  \rho^2= \rho^3, q^1=q^2+q^3\ ,$$

{\bf Case III:} $\hat Z^1 (1-\hat w^3)  \le \hat w^2- \hat w^3 $  and $\hat w^2 \ge \hat w^3$:
 $$q^2 = 0,\  q^3=q^1, \rho^1=\rho^3\ ,$$
 	
{\bf Case IV:} $\hat Z^1 (1-\hat w^2) \le  \hat w^3- \hat w^2$  and $\hat w^2 \le \hat w^3$:
  $$q^2=q^1,\ q^3 =0,\ \rho^1 = \rho^2\ .$$

One observes that the restrictions in the above cases in characteristic variables do lead to straightforward restrictions when writing them in macroscopic variables.
Additionally,  we have the physical constraints $0 \le q^1 \le \rho^1\le 1$. 
Writing the restrictions in macroscopic variables and using the above conditions in the respective cases, one obtains the following.
For  Case I one obtains with $\rho = \rho^1=\rho^2 = \rho^3$ that $$\hat Z^1 (1-\hat w^3)  \ge \hat w^2- \hat w^3$$ is equivalent to $$(q^2+q^3) (1-\rho +q^3 )\ge (q^3-q^2)(1-\rho+q^2+q^3)\ .$$
This gives 
$$
(q^2+q^3) (1-\rho +q^3 )\ge (q^3-q^2)(1-\rho+q^2+q^3)
$$
and
$$
q^2(2+q^3+q^2)\ge 0 \ ,
$$
which is obviously fulfilled.
 
Moreover 
$$\hat w^2 + \hat w ^3  \ge \hat Z^1$$
gives
$$ (2 \rho- q^1 )(1-w^1)\ge q^1$$
or
\begin{align*}
 (\rho^1+w^1)(1-w^1)
&= \rho^1 +w^1 -\rho^1 w^1 -(w^1)^2\\
&= 2 w^1 +q^1 -\rho^1 w^1 -(w^1)^2 
= w^1(2-\rho^1-w^1)+q^1
\ge q^1\ .
\end{align*}
This is again obvious from the physical bounds. 

For Case II with $\rho^1 = q^1, \rho^2= \rho^3=\rho$ we have the constraint
$$\hat w^2 + \hat w ^3  \le \hat Z^1\ .$$
This is equivalent to 
$$(2 \rho - q^1)\le \hat q^1$$
or 
$$\rho \le q^1\ .$$
Case III with $q^2 = 0, q^3=q^1, \rho^1=\rho^3=\rho$ has the constraints $\hat Z^1 (1-\hat w^3)  \le \hat w^2- \hat w^3 $   and $\hat w^2 \ge \hat w^3$. 
This is equivalent to 
$$
q^1 (1-\rho+q^1)  \le  (1-\rho+q^1)(\rho^2- \rho+q^1), \; \;  \rho^2 \ge  \rho -q^1 
$$
or
$$q^1   \le  \rho^2- \rho+q^1, \; \; \rho^2 -  \rho +q^1 \ge 0\ ,$$
which gives $$ \rho  \le  \rho^2\ .$$
Case IV is symmetric to case III.

Thus,  the above 4 cases  can be rewritten using  macroscopic variables:

{\bf Case I:} 
 $$
  \rho^1 = \rho^2= \rho^3 \ge q^1\  ,
 $$
 
{\bf Case II:} 
  $$
   \rho^2 = \rho^3 \le  q^1 = \rho^1\  ,
  $$
  
{\bf Case III:}
       $$q^1 \le  \rho^1 = \rho^3 \le \rho^2\ ,\quad   q^2 = 0\ ,$$
   
{\bf Case IV:}  	
      	  $$q^1 \le \rho^1 = \rho^2 \le \rho^3\ ,\quad q^3 =0 .$$

In  all cases, we have additionally the balance of fluxes $q^1=q^2+q^3$.

\begin{remark}
\label{rem1}
We note that the expression for $\rho^1$ can be rewritten as 
$$
\rho^1 = \max(q^1, \min(\rho^2,\rho^3))\ .
$$
\end{remark}

\subsection{Junction with diverging lanes and equal driver preferences}
\label{equalpref}
For this configuration several coupling conditions have been proposed for the macroscopic conservation law, see e.g. \cite{CGP05,L,Dag2,HK}.
These conditions rely on a given preference of the drivers, i.e. it is known a priory what percentage of the arriving cars will take road $2$ and which ones road $3$.
We consider a simplified situation with an equal distribution of the percentage of cars which have a preference for road $2$ or  $3$ respectively. 
This fixes already two equations at the junction
\begin{align}\label{eq:FIFO_kin_eq1}
f_1^2 = \frac{ f_1^1}{2} \qquad \text{ and }\qquad f_1^3 = \frac{f_1^1}{2}\ .
\end{align}
The discussion of the free space on road 1  provides again the missing  information. 
For example, in case the available space on road 2 is larger than on road 3, the free space available for the drivers is the one  on road 3 plus the space occupied by the cars driving towards road 2. 

If $1-\hat w^2 \ge 1-w^3$ or $w^2\le w^3$ or $\rho^2 \le \rho^3$:
\begin{align*}
 1-w^1= 1- f_0^1 = 1-f_0^3+ q^2= 1-\hat w^3 +q^2 =  1-\hat w^3 +\frac{q^1}{2}\ .
\end{align*}

If $1-\hat w^2 \le 1-w^3$ or $w^2\ge w^3$ or $\rho^2 \ge \rho^3$:
\begin{align}\label{eq:FIFO_kin_eq3}
1-w^1= 1- f_0^1 = 1-f_0^2+ q^3= 1-\hat w^2 +q^3 =  1-\hat w^2 +\frac{q^1}{2}\ .
\end{align}

These three equations (\ref{eq:FIFO_kin_eq1} - \ref{eq:FIFO_kin_eq3}) form the coupling conditions, if they yield solutions within the physcial bounds.
Note that the conservation of mass is guaranteed by \eqref{eq:FIFO_kin_eq1}.
	
We consider them again in characteristic variables. 
This leads for $w^2 \le w^3$ to 	
\begin{align*}
1-w^1 &= 1-\hat w^3+ \frac{\hat Z^1}{2} (1-w^1)\\
(1- \hat w^2) Z^2 &=  (1- w^1) \frac{\hat Z^1}{2}\\
(1- \hat w^3)  Z^3 &=  (1- w^1) \frac{\hat Z^1}{2}
\end{align*}
or 
\begin{align*}
  1-w^1=& \frac{1-\hat w^3}{1-\frac{\hat Z^1}{2} }\ ,\\
Z^2 =   \frac{1-\hat w^3}{1-\hat w^2}\frac{\frac{\hat Z^1}{2}}{1-\frac{\hat Z^1}{2}}
\ ,
\qquad &\qquad 
Z^3 =  \frac{\frac{\hat Z^1}{2}}{1-\frac{\hat Z^1}{2}}\ .
\end{align*}
This is well defined as long as 
$$
1-\hat w^3 \le 1-\frac{\hat Z^1}{2} \; \; \mbox{and}\; \; \frac{\hat Z^1}{2} \le 1-\hat w^2\ 
$$
or
$$
\hat w^3 \ge \frac{\hat Z^1}{2} \; \; \mbox{and}\; \; \hat w^2 \le 1 -\frac{\hat Z^1}{2}\ .
$$
For $w^2 \ge w^3$ one has symmetrically
\begin{align*}
1-w^1=& \frac{1-\hat w^2}{1-\frac{\hat Z^1}{2} }\ ,\\
Z^2 =  \frac{\frac{\hat Z^1}{2}}{1-\frac{\hat Z^1}{2}}
\ ,\qquad & \qquad
Z^3 =   \frac{1-\hat w^2}{1-\hat w^3}\frac{\frac{\hat Z^1}{2}}{1-\frac{\hat Z^1}{2}}
\end{align*}
These expressions are  well defined as long as 
$$
1-\hat w^2 \le 1-\frac{\hat Z^1}{2} \; \; \mbox{and}\; \; \frac{\hat Z^1}{2} \le 1-\hat w^3\ 
$$
or
$$
\hat w^2 \le \frac{\hat Z^1}{2} \; \; \mbox{and}\; \; \hat w^3 \le 1 -\frac{\hat Z^1}{2}\ .
$$
To complete the coupling conditions we have to consider the remaining cases $0\le w^2,w^3 \le \frac{\hat Z^1}{2}$ and $1-\frac{\hat Z^1}{2}\le w^2,w^3 \le 1$ and truncate the coupling conditions in a suitable way.
In  the first case we use $1-w^1 =1$. In the second case the above expression for $1-w^1$ can be used.
This leads to
\begin{align*}
1-w^1 =
\begin{cases}
\frac{1-\hat w^3}{1- \frac{\hat Z^1}{2}}&, \hat w^3 \ge \frac{\hat Z^1}{2},\hat w^2 \le \hat w^3\\
\frac{1-\hat w^2}{1- \frac{\hat Z^1}{2}}&, \hat w^2 \ge \frac{\hat Z^1}{2},\hat w^2 \ge \hat w^3\\
1 &, \hat w^2\le \frac{\hat Z^1}{2}, \hat w^3 \le \frac{\hat Z^1}{2}.
\end{cases}
\end{align*}
and
the associated $Z^2,Z^3$ as
$$
\hat Z^2 = \frac{1-w^1}{1-\hat w^2} \frac{\hat Z^1}{2}\ , \qquad \hat Z^3 = \frac{1-w^1}{1-\hat w^3} \frac{\hat Z^1}{2}\ .
$$

\begin{figure}[h!]
	\begin{center}
		\externaltikz{figurepref}{
	 	\begin{tikzpicture}[scale = 4]
	 	\def\dr{0.2}
	 	\node[below] at (0,0) {$(0,0)$};
	 	\node[below] at (1,0) {$(1,0)$};
	 	\node[right] at (1,1) {$(1,1)$};
	 	\node[left] at (0,1) {$(0,1)$};
	 		 		 	\node[below] at (0.15,0.2) {$\bf I$};
	 		 		 	\node[right] at (0.7,0.2) {$\bf II$};
	 		 		 	\node[left] at (0.4,0.8) {$\bf III$};
	 	\node[below] at (0.3,0) {$\frac{\hat Z^1}{2}$};
	 	\node[below] at (0.7,0) {$1-\frac{\hat Z^1}{2}$};
	 	\node[left] at (0,0.3) {$\frac{\hat Z^1}{2}$};
	 		 	\node[left] at (0,0.7) {$1-\frac{\hat Z^1}{2}$};
	 	\draw[line width=1pt](0.3,0)--(0.3,0.3);
	 	\draw[line width=1pt](0,0.3)--(0.3,0.3);
	 	\draw[->](0,0)--(1.2,0) node[below]{$\hat w^2$};
	 	\draw[->](0,0)--(0,1.2) node[left]{$\hat w^3$};
	 	\draw(1,0)--(1,1);
	 	\draw(0,1)--(1,1);
	 	\draw[line width=1pt](0.3,0.3)--(1,1);
	 	\end{tikzpicture}
	 	 }
	\end{center}
	\caption{Domain of definition for coupling conditions with equal drivers preference.}
	\label{fig:pref}
\end{figure}

In macroscopic variables  $\rho$ and $q$ we obtain from the definition of $1-w^1$  the following three cases, see Figure \ref{fig:pref}.

{\bf Case I:}  $\hat w^2 \le \frac{\hat Z^1}{2}$ and $ \hat w^3 \le \frac{\hat Z^1}{2}$:\qquad $\rho^1=q^1$,

{\bf Case II:}  $\hat w^2 \ge \frac{\hat Z^1}{2}$ and $ \hat w^2  \ge \hat w^3$:\qquad $\rho^1=\rho^3$,

{\bf Case III:}  $\hat w^3 \ge \frac{\hat Z^1}{2}$ and $\hat w^2  \le \hat w^3$:
\qquad $\rho^1=\rho^2$.\\

\noindent
A short computation shows that this is can be rewritten as

{\bf Case I:}  \vspace{-1.5em}
$$
\rho^2, \rho^3 \le q^1=\rho^1\ ,
$$

{\bf Case II:} \vspace{-1.5em}
$$
q^1  \le \rho^1 = \rho^3 \le \rho^2\ ,
$$

{\bf Case III:}  \vspace{-1.5em}
$$
q^1  \le \rho^1  = \rho^2\le \rho^3\ 
$$
together with the equations for the fluxes
$$
q^2= \frac{q^1}{2} = q^3.
$$

\begin{remark}
\label{rem2}
 This can be rewritten as
\begin{align*}
\rho^1 = \max(q^1, \max(\rho^2,\rho^3))\ ,
\end{align*}
compare Remark \ref{rem1} for the case without drivers preferences.
\end{remark}

\section{The layer equations and the half-Riemann problem for the conservation law}
\label{layerRP}
In this section we reconsider the kinetic layer equations and their asymptotic states and the half-Riemann problems for the conservation law, see \cite{BK20} for more details.

\subsection{Layer solutions for the kinetic  equations}
\label{kinlayer}

\subsubsection{Left layer}
Let the left boundary of the domain be located at $x=x_L$.
Starting from equation \eqref{macro0} and rescaling space as $y= \frac{x-x_L}{\epsilon}$ and neglecting higher order terms in $\epsilon$ we obtain the kinetic layer equations for the left boundary for $(\rho_L,q_L)$ and $y \in [0, \infty)$ as 
\begin{align*}
\begin{aligned}
\partial_y q_L &=0\\
\frac{q_L}{1-\rho_L}  \partial_y \rho_L  + (1-\frac{q_L}{1-\rho_L})\partial_y  q_L & =- \left(q_l-F(\rho_L) \right) \ .
\end{aligned}
\end{align*}
This yields 
\begin{align*}
\begin{aligned}
q_L =C\ ,\quad 
\partial_y \rho_L   = (1-\rho_L) \frac{F(\rho_L)-C}{C}\ .
\end{aligned}
\end{align*}
For $0<C < F(\rho^{\star})=\sigma$, where $\rho^{\star}$ denotes the point where the maximum of $F$ is attained.
The above problem has two relevant fix-points
$$\rho_{-} (C)\le \rho^{\star}\ ,\ 
\rho_+ (C)= \tau (\rho_-) \ge \rho^{\star}\ .$$
Here, $\tau(\rho)\neq \rho$ is defined by $F(\tau(\rho))= F(\rho)$.
The point $\rho_-$ is instable and $\rho_+ $ is stable. The domain of attraction of the stable fixpoint $\rho_+$  is   the interval $(\rho_-,1)$.

The third fixpoint $\rho=1$ is not relevant for the further matching procedure, since  in the macroscopic limit the maximal density requires $C=0$. 
In case $C=0$ we have the instable fixpoint $\rho_+ = 1$ and the stable fixpoint $\rho_- =0$ with domain of attraction $[0,1)$. 
Further, for $C=F(\rho^{\star})$ both fixpoints coincide, i.e. $\rho_- = \rho_+ = \rho^{\star}$, and all solutions with initial values above $\rho^{\star}$ converge towards $\rho^{\star}$, all other solutions diverge.

\begin{remark}
In case of the LWR model with $F(\rho) = \rho(1-\rho)$ we have with $C< \frac{1}{4}$ 
$$\rho_{\pm} (C) = \frac{1}{2} (1 \pm \sqrt{1-4 C})\ .$$
In case $C=\frac{1}{4}$ it holds $\rho_- = \rho_+ = \frac{1}{2}$. 
Moreover, $\tau(\rho) = 1-\rho$, as in Figure \ref{figfund}.
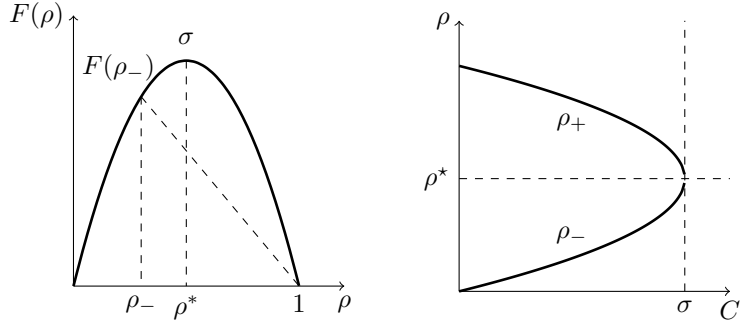
\begin{figure}[h]
	\center
	\externaltikz{LWR}{
			\begin{tikzpicture}[scale = 3]
			\def \rhobar {0.3}	
			\def \rhostar {0.5}
			\draw[->] (0,0)--(1.2,0) node[below]{$\rho$};
			\draw[->] (0,0)--(0,1.2) node[left]{$F(\rho)$}node at (0.2,0.97) {$F(\rho_-)$};
			\draw[black,line width=1pt,domain=0.0:1,smooth,variable=\x,] plot ({\x},{4*\x*(1-\x)}) ;
			\draw[dashed] (\rhobar,{4*\rhobar*(1-\rhobar)})--(\rhobar,0) node[below]{$\rho_-$};
			\draw[dashed] (\rhostar,{4*\rhostar*(1-\rhostar)})--(\rhostar,0) node[below]{$\rho^*$}node at (0.5,1.1) {$\sigma$};;
			\draw[dashed] (\rhobar,{4*\rhobar*(1-\rhobar)})--(1,0)
			node[below]{$1$};	;
			\end{tikzpicture}
			\hspace{0.5cm}
		\begin{tikzpicture}[scale = 3]	
			\draw[->] (0,0)--(1.2,0) node[below]{$C$};
			\draw[->] (0,0)--(0,1.2) node[left]{$\rho$};
			\draw[dashed] (0,0.5)--(1.2,0.5) ;
			\draw[dashed] (1.0,0)--(1.0,1.2) ;
			\node at (0.5,0.25) {$\rho_-$};
	  		\node at (0.5,0.75) {$\rho_+$};
		 	\node[left] at (0,0.5) {$\rho^\star$}; 
		 	\node[below] at (1.0,0) {$\sigma$}; 
			\draw[black,line width=1pt,domain=0.0:1.0,smooth,variable=\x,samples = 100] plot ({\x},{0.5*(1-sqrt(1-\x))});
			\draw[black,line width=1pt,domain=0.0:1.0,smooth,variable=\x,samples = 100] plot ({\x},{0.5*(1+sqrt(1-\x))});
		\end{tikzpicture}
	}
	\caption{Fundamental diagram, $F(\rho)$ and $\rho_\mp$.}
	\label{figfund}
\end{figure}
\end{remark}

\subsubsection{Right layer}
For the right boundary at $x_R$ a similar scaling $y=\frac{x_R-x}{\epsilon}$ gives the layer equations for $(\rho_R,q_R)$ and  $y\in[0, \infty)$ as
\begin{align*}
\begin{aligned}
q_R &=C\\
- \partial_y \rho_R   &= (1-\rho_R) \frac{F(\rho_R)-C}{C}\ .
\end{aligned}
\end{align*}
For  $0 < C < F(\rho^{\star})$
the above problem has again two  relevant fix points
$$\rho_{-}(C) \le \rho^{\star}\ ,\ 
\rho_+ (C) = \tau (\rho_-) \ge \rho^{\star}\; .$$
In this case 
$\rho_- $ is stable, $\rho_+ $ is instable. 
The domain of attraction of the  stable fixpoint $\rho_-$ is $[0,\rho_+)$.

For $C=F(\rho^{\star})=\sigma$ we have $\rho_-= \rho_+ = \rho^{\star}$ and 
all solutions with initial values below  $\rho^{\star}$ converge towards $\rho^{\star}$, all other solutions converge to not admissible states.
For $C=0$ the fixpoint $\rho_+ = 1$ is instable and $\rho_- =0$ is a stable fixpoint with domain of attraction $[0,1)$.

\subsubsection{Summary}
\label{summary}
In summary we have the following cases denoting with  $U$ the unstable fixpoints and with $S$ the stable ones.
Moreover, we use the notation $\rho_K$ for the  values $\rho^\infty_L$ and $\rho^\infty_R$ at infinity of the respective  layers 
and the notation $\rho_0$ for the respective values at $y=0$, i.e. $\rho_L(0)$ and $\rho_R(0)$.
\\
\paragraph{Layer Problem at the left boundary} 
\begin{align*}
&\left.\begin{array}{lll}
\rho_K = \rho_-(C) \quad &\Rightarrow\quad \rho_0 =\rho_-(C), &0 \le   C< \sigma
\end{array}\right\}
&\quad \text{(U)}\\
&\left.\begin{array}{lll}
\rho_K = \rho_+(C) \quad &\Rightarrow\quad \rho_0 \in (\rho_-(C),1), &0< C < \sigma\\
\rho_K = \rho^\star \quad &\Rightarrow\quad  \rho_0 \in [\rho^\star,1),& C=\sigma\\
\rho_K = 1 \quad &\Rightarrow\quad  \rho_0 \in (0,1],& C=0
\end{array}\right\}
&\quad \text{(S)}
\end{align*}

\paragraph{The Layer Problem at the right boundary} 			
\begin{align*}
&\left.\begin{array}{lll}
\rho_K = \rho_+(C) \quad &\Rightarrow\quad \rho_0 =\rho_+(C), &0 \le  C< \sigma 
\end{array}\right\}
&\quad \text{(U)}\\
&\left.\begin{array}{lll}
\rho_K = \rho_-(C) \quad &\Rightarrow\quad \rho_0 \in [0,\rho_+(C)), &0 < C < \sigma\\
\rho_K = \rho^\star \quad &\Rightarrow\quad  \rho_0 \in [0,\rho^\star],& C=\sigma\\
\rho_K = 0 \quad &\Rightarrow\quad  \rho_0 \in [0,1),& C=0
\end{array}\right\}
&\quad \text{(S)}
\end{align*}

In the following we use for the three cases  of the stable fixpoint (S) the notation
$$
\rho_K = \rho_+(C) \quad \Rightarrow\quad  \rho(0) \in \lceil \rho_-(C),1 \rfloor, 0 \le  C \le  \sigma
$$
for the left boundary and for the right boundary
$$
\rho_K = \rho_-(C) \quad \Rightarrow\quad  \rho(0) \in \lceil 0, \rho_+(C) \rfloor, 0 \le  C \le  \sigma\ .
$$

\subsection{Half-Riemann problems for the limit conservation law}
\label{Riemann}
Assuming the conditions above on $F$ the solution to a Riemann problem of the limit conservation law  $ \partial_t \rho + \partial_x F(\rho)=0$ is easily obtained.
The possible states $\rho_K$ for a given value $\rho_B $ of a half-Riemann problem with ingoing waves (shocks and rarefaction waves) at left and right boundary are summarized in the following:\\

\paragraph{The half-Riemann Problem at the left boundary}
\begin{align*}
\rho_B&\leq \rho^\star \ (\text{RP 1})  \quad &\Rightarrow\quad \rho_K &\in [0,\rho^\star ]\\
\rho_B&> \rho^\star \ (\text{RP 2})  \quad &\Rightarrow\quad \rho_K &\in [0,\tau(\rho_B)]\cup\{\rho_B\}
\end{align*}		

\paragraph{The half-Riemann Problem at the right boundary}		
\begin{align*}
\rho_B&\geq  \rho^\star \ (\text{RP 1})  \quad &\Rightarrow\quad \rho_K &\in [ \rho^\star ,1]\\
\rho_B&<  \rho^\star  \ (\text{RP 2})  \quad &\Rightarrow\quad \rho_K &\in \{\rho_B\}\cup [\tau(\rho_B),1]
\end{align*}

These set will allow waves to emerge from the junction into the domains.

\section{Macroscopic coupling conditions: diverging lanes with no driver preferences}
\label{macroscopiccc}

For the determination of the coupling conditions for the macroscopic equations we  investigate first  the kinetic layers at the nodes coupled to each other via the coupling conditions and determine their asymptotic states. Then, we   match these  results to Riemann solutions of the  macroscopic problems on each of the roads.

Assuming   the  initial states $\rho_B^1,\rho_B^2,\rho_B^3$  on all three roads to be given, we have to determine the new states $\rho_K^1,\rho_K^2$ and $\rho_K^3$ at the node. On the one hand $\rho_K^1,\rho_K^2$ and $\rho_K^3$ are the asymptotic states of the respective layer problems, on the other hand they are the right (for road 1 and 2) or left (for road 3) states of the half-Riemann problems with $\rho_B^1,\rho_B^2,\rho_B^3$
as the corresponding left (road 1 and 2) or right  state (road 3).
We have to consider eight different configurations of Riemann problems.
For  each of them all possible combinations with stable or unstable layer solutions have to be discussed.
Not admissible combinations are not listed.
The proof of the following theorem is given in section \ref{technical}.
For the discussion of the coupling of the layer solutions we refer to \ref{layerproof} and for the matching of layer solutions and half Riemann problems to  
\ref{proof}.
One obtains

\begin{theorem}
	\label{theoremdiverge}
Starting from the kinetic coupling conditions for drivers without preferences  in subsection \ref{nopref} the asymptotic derivation of the  coupling conditions  for the macroscopic equations gives  the following  cases using the notation 
RP1/2-1/2-1/2 for the respective combination of the half Riemann problems.

\noindent{\bf Case 1, RP1-1-1} $\rho_B^1 \ge \rho^\star\; ,\; \rho_B^2 \le \rho^\star\; ,\; \rho_B^3 \le  \rho^\star $.

One obtains $$C^1 = \sigma\; , \; C^2=C^3=\frac{\sigma}{2}\ .$$

\noindent{\bf Case 2, RP1-1-2} $\rho_B^1 \ge \rho^\star , \rho_B^2 \le \rho^\star , \rho_B^3 \ge  \rho^\star $.

This gives for 
\begin{align*}
2 F(\rho_B^3)) \le \sigma&:C^1 =  \sigma\;,\; C^3 = F(\rho_B^3)\\
 2 F(\rho_B^3)) \ge \sigma &:C^1 = \sigma \;,\; C^2=C^3= \frac{\sigma}{2}\ .
\end{align*}

\noindent{\bf Case 3, RP1-2-1} $\rho_B^1 \ge \rho^\star \;,\; \rho_B^2 \ge \rho^\star \;,\; \rho_B^3 \le  \rho^\star $. 

This gives  for 
\begin{align*}
 2 F(\rho_B^2) \le \sigma&:C^1 = \sigma \;,\; C^2=C^3 = \frac{\sigma}{2}\\
  2 F(\rho_B^2) \ge \sigma&:C^1 =  \sigma\;,\; C^2 = F(\rho_B^2)\ . 
\end{align*}

\noindent{\bf Case 4, RP2-1-1} $\rho_B^1 \le \rho^\star\; ,\; \rho_B^2 \le \rho^\star \;,\; \rho_B^3 \le  \rho^\star $. 

This gives $$C^1 = F(\rho_B^1) \ \mbox{and} \ C^2= C^3 =\frac{ F(\rho_B^1)}{2}\ .$$

\noindent{\bf Case 5, RP1-2-2} $\rho_B^1 \ge \rho^\star \;,\; \rho_B^2 \ge \rho^\star \;,\; \rho_B^3 \ge  \rho^\star $. 

This gives for 
\begin{align*}
F(\rho_B^2)+  F(\rho_B^3) \le \sigma&: 
C^2 = F(\rho_B^2) , C^3 = F(\rho_B^3)\ , \\
F(\rho_B^2) \ge \frac{\sigma}{2}\;,\; F(\rho_B^3) \ge \frac{\sigma}{2} &: 	
C^1 = \sigma\;,\;C^2= \frac{\sigma}{2},C^3= \frac{\sigma}{2}\ ,\\
F(\rho_B^2)+  F(\rho_B^3) \ge \sigma,F(\rho_B^2) \ge \frac{\sigma}{2}&: 	C^1 = \sigma \;,\; C^2  = \sigma- F(\rho_B^3), C^3 = F(\rho_B^3)\ ,\\ 
F(\rho_B^2)+  F(\rho_B^3) \ge \sigma,F(\rho_B^3) \ge \frac{\sigma}{2}&: 	C^1  = \sigma\;,\; C^2 = F(\rho_B^2),C^3 =\sigma- F(\rho_B^2)\ .
\end{align*}	

\noindent{\bf Case 6, RP2-1-2} $\rho_B^1 \le \rho^\star \;,\; \rho_B^2 \le \rho^\star \;,\; \rho_B^3 \ge  \rho^\star $. 

This gives  for 
\begin{align*}
 F(\rho_B^1)  \le 2  F(\rho_B^3)&:
C^1 =  F(\rho_B^1)\;,\;C^2=C^3= \frac{C^1}{2} = \frac{F(\rho_B^1)}{2}\ ,\\
2 F(\rho_B^3)\le   F(\rho_B^1)&: C^1 =  F(\rho_B^1)\;,\; C^3 = F(\rho_B^3)\;,\; C^2=F(\rho_B^1)- F(\rho_B^3)\ .
\end{align*}

\noindent{\bf Case 7, RP2-2-1} $\rho_B^1 \le \rho^\star \;,\; \rho_B^2 \ge \rho^\star \;,\; \rho_B^3 \le  \rho^\star $. 

This gives for 
\begin{align*}
 F(\rho_B^1)  \le 2  F(\rho_B^2)&:
C^1 =  F(\rho_B^1)\;,\;C^2=C^3= \frac{C^1}{2} = \frac{F(\rho_B^1)}{2}\ ,\\
  2 F(\rho_B^2)\le   F(\rho_B^1)&:C^1 =  F(\rho_B^1)\;,\; C^2 = F(\rho_B^2)\;,\; C^3=F(\rho_B^1)- F(\rho_B^2)\ .
  \end{align*}

\noindent{\bf Case 8, RP2-2-2} $\rho_B^1 \le \rho^\star \;,\; \rho_B^2 \ge \rho^\star \;,\; \rho_B^3 \ge  \rho^\star $.

This gives for 
\begin{align*}
 F(\rho_B^2) +   F(\rho_B^3)\le F(\rho_B^1)&:
		C^1 = F(\rho_B^2)+ F(\rho_B^3) \;,\; C^2 =F(\rho_B^2)\;,\; C^3 = F(\rho_B^3)\ ,\\
		 \frac{F(\rho_B^1)}{2}\le \min ( F(\rho_B^2)\;,\;  F(\rho_B^3))&:
C^2= C^3= \frac{F(\rho_B^1)}{2}\;,\; C^1 = F(\rho_B^1)\ ,\\
F(\rho_B^2)+ F(\rho_B^3) \ge F(\rho_B^1)&: \\
2 F(\rho_B^3)\le   F(\rho_B^1)&:
C^1 =  F(\rho_B^1)\;,\; C^2= F(\rho_B^1) -F(\rho_B^3),C^3 = F(\rho_B^3)\ ,\\
  F(\rho_B^2)\le 2  F(\rho_B^1)&:
C^1  = F(\rho_B^1)\;,\; C^2 =  F(\rho_B^2) \;,\;C^3  =F(\rho_B^1)- F(\rho_B^2)\ .
\end{align*}

\end{theorem}

\subsection{Supply-Demand formulation of the coupling conditions}

We use the supply-demand representation and denote the sets of valid resulting fluxes $C^i$ by $\Omega^i$, compare \cite{CGP05,L,Dag1,Dag2,HR}. 

The
sets $\Omega^i$ are
for the incoming road $i=1,2$
\begin{align*}
\rho_B^i \le \rho^\star \Rightarrow \Omega^i = [0, F(\rho_B^i)] &&\text{ and }&&
\rho_B^i \ge \rho^\star \Rightarrow \Omega^i = [0, \sigma]\ .	 
\end{align*}
For the outgoing road $i=3$
\begin{align*}
\rho_B^i \le\rho^\star \Rightarrow \Omega^i = [0, \sigma]  &&\text{ and }&&
\rho_B^i\ge  \rho^\star \Rightarrow \Omega^i = [0, F(\rho_B^i)]\ .	 
\end{align*}
We define 
$c^i$ such that $\Omega^i = [0, c^i]$.
Rewriting the above conditions using this notation gives the following.

\noindent{\bf Case 1, RP1-1-1.} 
This is a case with 
$C^1 =c^1, C^2 = C^3= \frac{c^1}{2}$.

\noindent{\bf Case 2, RP1-1-2.}  
We have   two cases. 
 $ C^1= c^1, C^2 =C^3 =  \frac{c^1}{2}$, if $c^3 \ge \frac{c^1}{2}$ and  $ C^1= c^1, C^2 = c^1- c^3, C^3 =  c^3$, if $c^3 \le \frac{c^1}{2}$.

\noindent{\bf Case 3, RP1-2-1} 
We have    two cases.  $ C^1= c^1, C^2 = c^2, C^3 =  c^1- c^2$, if $c^2 \le \frac{c^1}{2}$ and  
 $ C^1=c^1, C^2 = C^3 =  \frac{c^1}{2}$, if $c^2 \ge \frac{c^1}{2}$.

\noindent{\bf Case 4,	RP2-1-1} 
This is a case with 
$
C^1=c^1\;,\; C^2=C^3 =  \frac{c^1}{2}$.

\noindent{\bf Case 5, RP1-2-2} 
We have four  cases:
\begin{align*}
c^2+c^3 \le c^1&:  C^2=c^2\;,\; C^3 =c^3\;,\; C^1= c^2+ c^3\ ,\\
c^2 +c^3 \ge c^1, c^2 \ge \frac{c^1}{2}, c^2\ge \frac{c^1}{2}&:  C^1= c^1,C^2 =C^3=\frac{c^1}{2}\ ,\\
c^2 +c^3 \ge c^1\;,\; c^3,c^2 \ge \frac{c^1}{2}\;,\; c^3\le \frac{c^1}{2}&:  C^1= c^1\;,\; C^2 = c^1-c^3\;,\;C^3 =c^3\ ,\\
c^2+c^3 \ge c^1\;,\; c^3\;,\;c^2 \le \frac{c^1}{2}\;,\; c^2\ge \frac{c^1}{2}&:  C^1= c^1\;,\; C^2 =c^2\;,\; C^3=c^1-c^2\ .
\end{align*}

\noindent{\bf Case 6, RP2-1-2}
This case is the same as Case 2.

\noindent{\bf Case 7,	RP2-2-1}  
This is the same as Case 3.

\noindent{\bf Case 8, RP2-2-2}
This is again Case 5.

Summarizing this leads to only 4 different cases:

{\bf Case A:} \vspace{-1em}
\begin{align*}
c^2+c^3 \le c^1&:  C^1= c^2+c^3\;,\; C^2 =c^2\;,\; C^3=c^3\ ,
\end{align*}

{\bf Case B:}\vspace{-1em}
\begin{align*}
c^2 +c^3 \ge c^1\;,\;c^1 \le 2 c^2\;,\; c^1 \le 2 c^3&:  C^1= c^1\;,\; C^2 =C^3= \frac{c^1}{2}\ ,
\end{align*}

{\bf Case C:}\vspace{-1em}
\begin{align*}
c^2 +c^3 \ge c^3\;,\;c^1 \ge 2c^2\;,\; c^1\le 2 c^3&:  C^1= c^1\;,\; C^2 = c^2\;,\; C^3=c^1-c^2\ ,
\end{align*}

{\bf Case D:}\vspace{-1em}
\begin{align*}
c^2+c^3 \ge c^1\;,\;c^1 \le 2 c^2\;,\;c^1\ge 2 c^3&:  C^1= c^1\;,\; C^2 =c^1-c^3\;,\; C^3 =c^3\ .
\end{align*}

The cases B)-D) can be rewritten as
\begin{align*}
C^1= c^1\;,\; C^i = \min\left(c^i,c^1-\min\left(c^2,c^3,\frac{c^1}{2}\right)\right)
\qquad i = 2,3\ .
\end{align*}

This yields the limit coupling conditions for the conservation law without drivers preferences at the node.

\subsection{Macroscopic coupling conditions: diverging lanes with equal \\driver preferences}
\label{Appendix:FIFO versus  NONFIFO models}

We use the same  notation as in the previous section, i.e. we
define $C_{i},c_i$ and the sets $\Omega_i=[0, c_i]$ as above, depending on
whether incoming or outgoing roads are considered.
The kinetic  conditions for drivers with equal preferences for each of the two outgoing lanes from Section \ref{equalpref} lead 
in the limit to the macroscopic coupling conditions
\begin{align*}
 C^1 = \min\Big(  c^1, 2  c^2 , 2c^3 \Big)
\end{align*}
and
$$
C^2 = C^3= \frac{C^1}{2}\ .
$$
We refer to \cite{CGP05,Dag1,Dag2} for scalar traffic models   on networks  with similiar coupling  conditions.
For other conditions
treating  situations with drivers preferences, we refer to  \cite{L,HK}.
Such  models can be  derived from suitable kinetic coupling conditions in a similiar way.

\section{Numerical results}
\label{Numerical results}

In this section we show some numerical examples to confirm the analytically derived coupling conditions. 
In the examples we compare the numerical solution of the kinetic model \eqref{eq:lindeg+relax} to the solution of the LWR model with the respective coupling conditions.
Both equations are approximated with a Godunov scheme and each edge is discretized with $1000$ cells. 
All solutions are shown at $T=0.9$.
The relaxation in the discrete velocity model is $\varepsilon = 0.001$.
As initial conditions the densities are chosen constant on the edges for both models equally, $Z$ in the kinetic model is set in equilibrium with the source term. 

In the figures the solution on the edges are shown on the left, on the right hand side a zoom close to the junction displays the possible layers. 

\subsection{Diverging lanes without driver preferences}
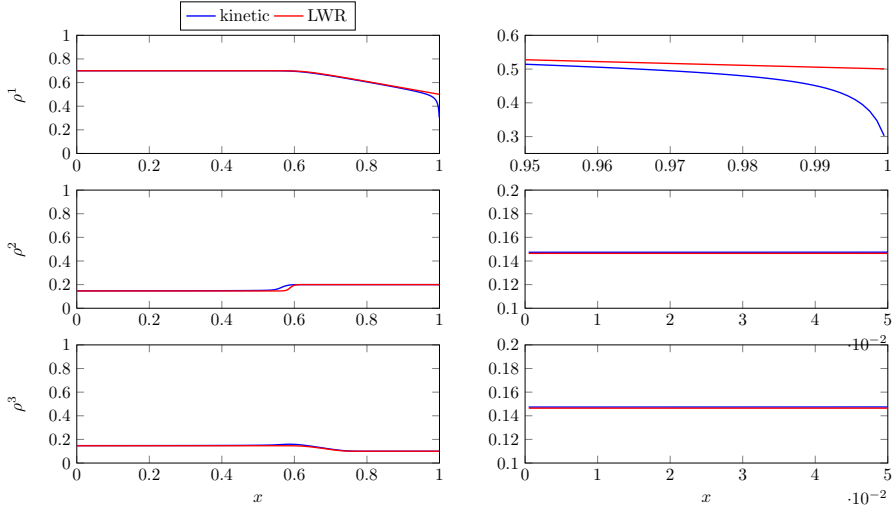
\begin{figure}[h!]
	\externaltikz{split_case4}{
		\begin{tikzpicture}[scale=0.65]
		\begin{groupplot}[
				group style={group size=2 by 3, vertical sep = 0.75cm, horizontal sep = 1.75cm},
				width = 9cm,
				height = 4cm,
				xmin = -0.0, xmax = 1.0,
				ymin = 0.0, ymax = 1.0,
				legend style = {at={(0.5,1)},xshift=0.2cm,yshift=0.1cm,anchor=south},
				legend columns= 3,			
				]

		\nextgroupplot[ ylabel = $\rho^1$]
			\addplot[color = blue,thick] file{Data/split_Lindeg_rho_1ex4_eps0001.txt};
			\addlegendentry{kinetic}
			\addplot[color = red,thick] file{Data/split_LWR_rho_1ex4.txt};
			\addlegendentry{LWR}
		\nextgroupplot[
				xmin = 0.95, xmax = 1.0,
				ymin = 0.25, ymax = 0.6]	
			\addplot[color = blue,thick] file{Data/split_Lindeg_rho_1ex4_eps0001.txt};
			\addplot[color = red,thick] file{Data/split_LWR_rho_1ex4.txt};
		
		\nextgroupplot[ ylabel = $\rho^2$]
			\addplot[color = blue,thick] file{Data/split_Lindeg_rho_2ex4_eps0001.txt};
			\addplot[color = red,thick] file{Data/split_LWR_rho_2ex4.txt};
		\nextgroupplot[
				xmin = 0.0, xmax = 0.05,
				ymin = 0.1, ymax = 0.2]	
			\addplot[color = blue,thick] file{Data/split_Lindeg_rho_2ex4_eps0001.txt};
			\addplot[color = red,thick] file{Data/split_LWR_rho_2ex4.txt};
		
		\nextgroupplot[ xlabel =  $x$, ylabel = $\rho^3$]
			\addplot[color = blue,thick] file{Data/split_Lindeg_rho_3ex4_eps0001.txt};
			\addplot[color = red,thick] file{Data/split_LWR_rho_3ex4.txt};
		\nextgroupplot[ xlabel =  $x$, 
				xmin = 0.0, xmax = 0.05,
				ymin = 0.1, ymax = 0.2]	
			\addplot[color = blue,thick] file{Data/split_Lindeg_rho_3ex4_eps0001.txt};
			\addplot[color = red,thick] file{Data/split_LWR_rho_3ex4.txt};
		\end{groupplot}
		
		\end{tikzpicture}
	}
	\caption{Diverging lanes example 1: $\rho^1 = 0.7$, $\rho^2 = 0.2$, $\rho^3 = 0.1$ }
	\label{fig:Split_case4}
\end{figure}

For the junction without driver preferences	we consider the initial conditions $\rho^1 = 0.7$, $\rho^2 = 0.2$ and $\rho^3 = 0.1$.
We are in the situation of Case 1 with $\rho_0^1=0.25, \rho_0^2=\rho_0^3=0.1464$.
In the incoming road the maximal flux at the junction generates a rarefaction wave, see Figure \ref{fig:Split_case4}.
This flux is distributed onto the outgoing roads, such that a small shock and a small rarefaction wave arise. 
On the right hand side we observe that a  layer forms in the first road but not in the two exiting ones.   
We observe that the equality  $\rho^2=\rho^3$ is still valid  on the macroscopic level.
Furthermore we note that in the kinetic model the shock on road $2$ is slightly behind the macroscopic one. 
It has the same speed as the macroscopic shock, but is  slightly delayed, since the kinetic model needs a few time steps to establish the correct states at the junction. 
Such initial layer problems decrease with decreasing $\varepsilon$ and increasing numerical resolution.

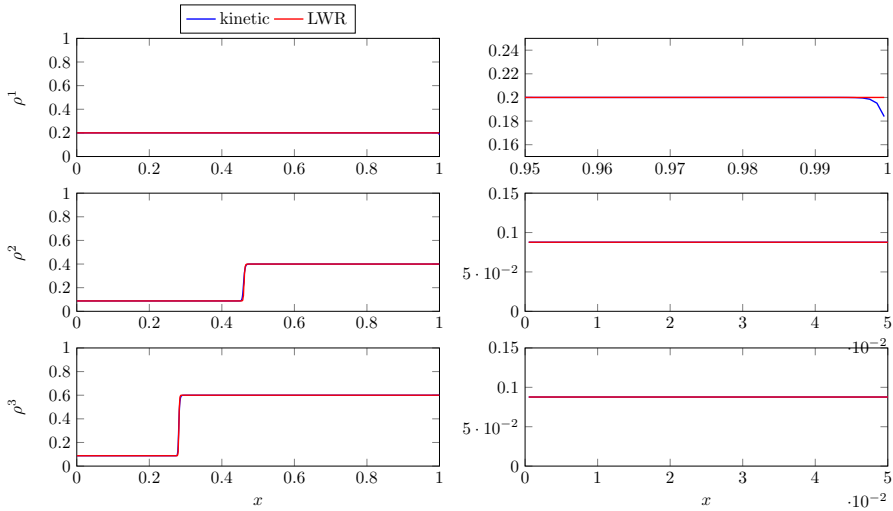
\begin{figure}[h!]
	\externaltikz{split_case5}{
		\begin{tikzpicture}[scale=0.65]
		\begin{groupplot}[
		group style={group size=2 by 3, vertical sep = 0.75cm, horizontal sep = 1.75cm},
		width = 9cm,
		height = 4cm,
		xmin = -0.0, xmax = 1.0,
		ymin = 0.0, ymax = 1.0,
		legend style = {at={(0.5,1)},xshift=0.2cm,yshift=0.1cm,anchor=south},
		legend columns= 3,			
		]
		
		\nextgroupplot[ ylabel = $\rho^1$]
		\addplot[color = blue,thick] file{Data/split_Lindeg_rho_1ex5_eps0001.txt};
		\addlegendentry{kinetic}
		\addplot[color = red,thick] file{Data/split_LWR_rho_1ex5.txt};
		\addlegendentry{LWR}
		\nextgroupplot[
		xmin = 0.95, xmax = 1.0,
		ymin = 0.15, ymax = 0.25]	
		\addplot[color = blue,thick] file{Data/split_Lindeg_rho_1ex5_eps0001.txt};
		\addplot[color = red,thick] file{Data/split_LWR_rho_1ex5.txt};
		
		\nextgroupplot[ ylabel = $\rho^2$]
		\addplot[color = blue,thick] file{Data/split_Lindeg_rho_2ex5_eps0001.txt};
		\addplot[color = red,thick] file{Data/split_LWR_rho_2ex5.txt};
		\nextgroupplot[
		xmin = 0.0, xmax = 0.05,
		ymin = 0.0, ymax = 0.15]	
		\addplot[color = blue,thick] file{Data/split_Lindeg_rho_2ex5_eps0001.txt};
		\addplot[color = red,thick] file{Data/split_LWR_rho_2ex5.txt};
		
		\nextgroupplot[ xlabel =  $x$, ylabel = $\rho^3$]
		\addplot[color = blue,thick] file{Data/split_Lindeg_rho_3ex5_eps0001.txt};
		\addplot[color = red,thick] file{Data/split_LWR_rho_3ex5.txt};
		\nextgroupplot[ xlabel =  $x$, 
		xmin = 0.0, xmax = 0.05,
		ymin = 0.0, ymax = 0.15]	
		\addplot[color = blue,thick] file{Data/split_Lindeg_rho_3ex5_eps0001.txt};
		\addplot[color = red,thick] file{Data/split_LWR_rho_3ex5.txt};
		\end{groupplot}
		
		\end{tikzpicture}
	}
	\caption{Diverging lanes example 2: $\rho^1 = 0.2$, $\rho^2 = 0.4$, $\rho^3 = 0.6$ }
	\label{fig:Split_case5}
\end{figure}

In the second example with the initial values $\rho^1 = 0.2$, $\rho^2 = 0.4$ and $\rho^3 = 0.6$ only few cars arrive at the junction.
As shown in Figure \ref{fig:Split_case5} these cars are distributed equally onto the outgoing roads.
Thus two shock waves form and move to the right.
A layer forms only on road $1$.
In this case we are in the situation of Case 6, subcase 1 with $\rho_0^1=F(\rho_B^1)= 0.16$ and $\rho_0^2=\rho_0^3=\rho_-(F(\rho_B^1/2)= 0.087$.

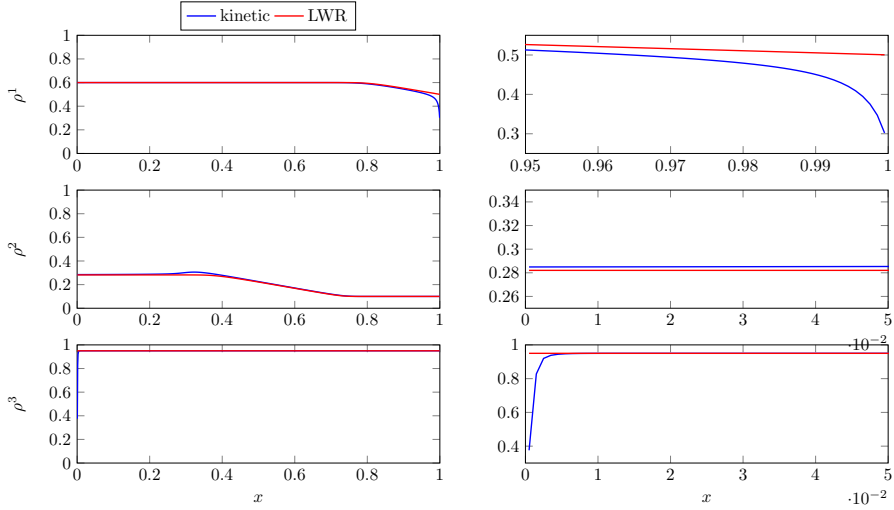
\begin{figure}[h!]
	\externaltikz{split_case6}{
		\begin{tikzpicture}[scale=0.65]
		\begin{groupplot}[
		group style={group size=2 by 3, vertical sep = 0.75cm, horizontal sep = 1.75cm},
		width = 9cm,
		height = 4cm,
		xmin = -0.0, xmax = 1.0,
		ymin = 0.0, ymax = 1.0,
		legend style = {at={(0.5,1)},xshift=0.2cm,yshift=0.1cm,anchor=south},
		legend columns= 3,			
		]
		
		\nextgroupplot[ ylabel = $\rho^1$]
		\addplot[color = blue,thick] file{Data/split_Lindeg_rho_1ex6_eps0001.txt};
		\addlegendentry{kinetic}
		\addplot[color = red,thick] file{Data/split_LWR_rho_1ex6.txt};
		\addlegendentry{LWR}
		\nextgroupplot[
		xmin = 0.95, xmax = 1.0,
		ymin = 0.25, ymax = 0.55]	
		\addplot[color = blue,thick] file{Data/split_Lindeg_rho_1ex6_eps0001.txt};
		\addplot[color = red,thick] file{Data/split_LWR_rho_1ex6.txt};
		
		\nextgroupplot[ ylabel = $\rho^2$]
		\addplot[color = blue,thick] file{Data/split_Lindeg_rho_2ex6_eps0001.txt};
		\addplot[color = red,thick] file{Data/split_LWR_rho_2ex6.txt};
		\nextgroupplot[
		xmin = 0.0, xmax = 0.05,
		ymin = 0.25, ymax = 0.35]	
		\addplot[color = blue,thick] file{Data/split_Lindeg_rho_2ex6_eps0001.txt};
		\addplot[color = red,thick] file{Data/split_LWR_rho_2ex6.txt};
		
		\nextgroupplot[ xlabel =  $x$, ylabel = $\rho^3$]
		\addplot[color = blue,thick] file{Data/split_Lindeg_rho_3ex6_eps0001.txt};
		\addplot[color = red,thick] file{Data/split_LWR_rho_3ex6.txt};
		\nextgroupplot[ xlabel =  $x$, 
		xmin = 0.0, xmax = 0.05,
		ymin = 0.3, ymax = 1]	
		\addplot[color = blue,thick] file{Data/split_Lindeg_rho_3ex6_eps0001.txt};
		\addplot[color = red,thick] file{Data/split_LWR_rho_3ex6.txt};
		\end{groupplot}
		
		\end{tikzpicture}
	}
	\caption{Diverging lanes example 3: $\rho^1 = 0.6$, $\rho^2 = 0.1$, $\rho^3 = 0.95$ }
	\label{fig:Split_case6}
\end{figure}

In Figure \ref{fig:Split_case6} the results with the initial conditions $\rho^1 = 0.6$, $\rho^2 = 0.1$ and $\rho^3 = 0.95$ are shown.
This is a sitiuation as in Case 2, subcase 1 with $\rho_0^1=\rho_0^2=\rho_0^3=\rho_-(\sigma-F(\rho_B^3))=0.2821$.
Since the traffic on road $3$ is dense only very few cars enter there. 
Most of the vehicles  enter into  road $2$.
In the kinetic solution we  observe two layers, one interacting with the rarefaction wave on road $1$ and one due to the ingoing characteristics on road $3$.
Note that the kinetic solution is very close to the macroscopic ones although the layer in road $3$ has to cover a range of more than $0.6$.
	
\subsection{Diverging lanes with driver preferences}
	
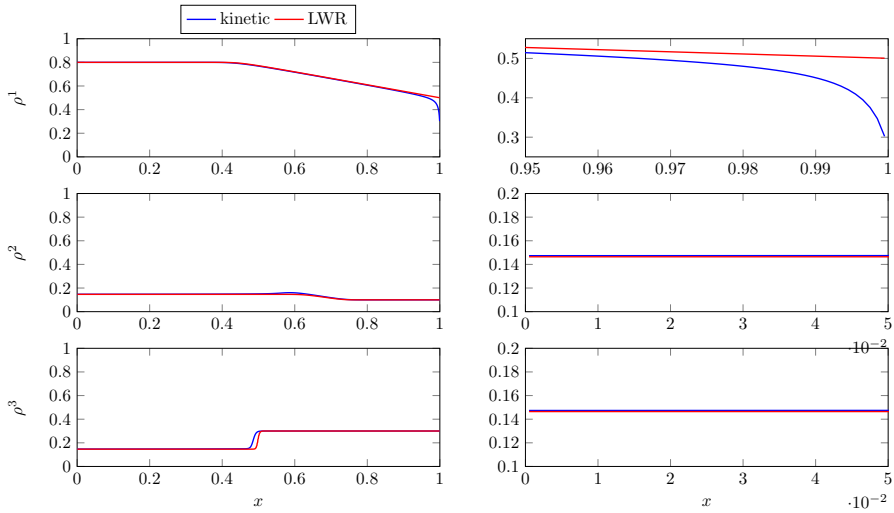
\begin{figure}[h!]
	\externaltikz{split_case7}{
		\begin{tikzpicture}[scale=0.65]
		\begin{groupplot}[
		group style={group size=2 by 3, vertical sep = 0.75cm, horizontal sep = 1.75cm},
		width = 9cm,
		height = 4cm,
		xmin = -0.0, xmax = 1.0,
		ymin = 0.0, ymax = 1.0,
		legend style = {at={(0.5,1)},xshift=0.2cm,yshift=0.1cm,anchor=south},
		legend columns= 3,			
		]
		
		\nextgroupplot[ ylabel = $\rho^1$]
		\addplot[color = blue,thick] file{Data/split_Lindeg_rho_1ex7_eps0001.txt};
		\addlegendentry{kinetic}
		\addplot[color = red,thick] file{Data/split_LWR_rho_1ex7.txt};
		\addlegendentry{LWR}
		\nextgroupplot[
		xmin = 0.95, xmax = 1.0,
		ymin = 0.25, ymax = 0.55]	
		\addplot[color = blue,thick] file{Data/split_Lindeg_rho_1ex7_eps0001.txt};
		\addplot[color = red,thick] file{Data/split_LWR_rho_1ex7.txt};
		
		\nextgroupplot[ ylabel = $\rho^2$]
		\addplot[color = blue,thick] file{Data/split_Lindeg_rho_2ex7_eps0001.txt};
		\addplot[color = red,thick] file{Data/split_LWR_rho_2ex7.txt};
		\nextgroupplot[
		xmin = 0.0, xmax = 0.05,
		ymin = 0.1, ymax = 0.2]	
		\addplot[color = blue,thick] file{Data/split_Lindeg_rho_2ex7_eps0001.txt};
		\addplot[color = red,thick] file{Data/split_LWR_rho_2ex7.txt};
		
		\nextgroupplot[ xlabel =  $x$, ylabel = $\rho^3$]
		\addplot[color = blue,thick] file{Data/split_Lindeg_rho_3ex7_eps0001.txt};
		\addplot[color = red,thick] file{Data/split_LWR_rho_3ex7.txt};
		\nextgroupplot[ xlabel =  $x$, 
		xmin = 0.0, xmax = 0.05,
		ymin = 0.1, ymax = 0.2]	
		\addplot[color = blue,thick] file{Data/split_Lindeg_rho_3ex7_eps0001.txt};
		\addplot[color = red,thick] file{Data/split_LWR_rho_3ex7.txt};
		\end{groupplot}
		
		\end{tikzpicture}
	}
	\caption{Diverging with driver preferences, example $1$: $\rho^1 = 0.8$, $\rho^2 = 0.1$, $\rho^3 = 0.3$ }
	\label{fig:Split_case7}
\end{figure}

For the junction with driver preferences we consider two slightly different examples. 	
First, with the initial conditions $\rho^1 = 0.8$, $\rho^2 = 0.1$ and $\rho^3 = 0.3$ there is enough space in both outgoing roads such that the maximal flow can be established, as shown in Figure \ref{fig:Split_case7}. We have $C^1=c^1$.
As the preferences for the two roads are equal, also the densities on road $2$ and $3$ are identical. 

\begin{figure}[h!]
	\externaltikz{split_case8}{
		\begin{tikzpicture}[scale=0.65]
		\begin{groupplot}[
		group style={group size=2 by 3, vertical sep = 0.75cm, horizontal sep = 1.75cm},
		width = 9cm,
		height = 4cm,
		xmin = -0.0, xmax = 1.0,
		ymin = 0.0, ymax = 1.0,
		legend style = {at={(0.5,1)},xshift=0.2cm,yshift=0.1cm,anchor=south},
		legend columns= 3,			
		]
		
		\nextgroupplot[ ylabel = $\rho^1$]
		\addplot[color = blue,thick] file{Data/split_Lindeg_rho_1ex8_eps0001.txt};
		\addlegendentry{kinetic}
		\addplot[color = red,thick] file{Data/split_LWR_rho_1ex8.txt};
		\addlegendentry{LWR}
		\nextgroupplot[
		xmin = 0.95, xmax = 1.0,
		ymin = 0.7, ymax = 0.8]	
		\addplot[color = blue,thick] file{Data/split_Lindeg_rho_1ex8_eps0001.txt};
		\addplot[color = red,thick] file{Data/split_LWR_rho_1ex8.txt};
		
		\nextgroupplot[ ylabel = $\rho^2$]
		\addplot[color = blue,thick] file{Data/split_Lindeg_rho_2ex8_eps0001.txt};
		\addplot[color = red,thick] file{Data/split_LWR_rho_2ex8.txt};
		\nextgroupplot[
		xmin = 0.0, xmax = 0.05,
		ymin = 0.75, ymax = 0.95]	
		\addplot[color = blue,thick] file{Data/split_Lindeg_rho_2ex8_eps0001.txt};
		\addplot[color = red,thick] file{Data/split_LWR_rho_2ex8.txt};
		
		\nextgroupplot[ xlabel =  $x$, ylabel = $\rho^3$]
		\addplot[color = blue,thick] file{Data/split_Lindeg_rho_3ex8_eps0001.txt};
		\addplot[color = red,thick] file{Data/split_LWR_rho_3ex8.txt};
		\nextgroupplot[ xlabel =  $x$, 
		xmin = 0.0, xmax = 0.05,
		ymin = 0.05, ymax = 0.15]	
		\addplot[color = blue,thick] file{Data/split_Lindeg_rho_3ex8_eps0001.txt};
		\addplot[color = red,thick] file{Data/split_LWR_rho_3ex8.txt};
		\end{groupplot}
		
		\end{tikzpicture}
	}
	\caption{Diverging with driver preferences, example $2$: $\rho^1 = 0.6$, $\rho^2 = 0.9$, $\rho^3 = 0.0$ }
	\label{fig:Split_case8}
\end{figure}
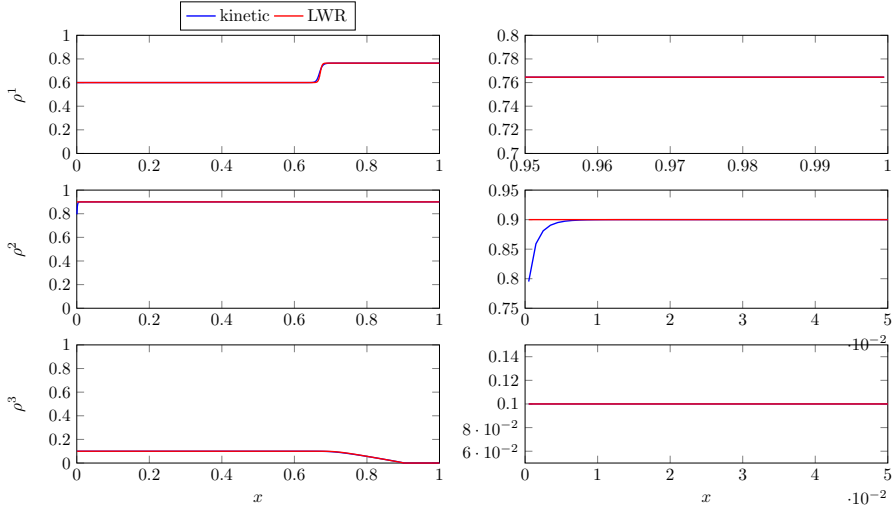

In Figure \ref{fig:Split_case8} the solutions corresponding to initial values $\rho^1 = 0.6$, $\rho^2 = 0.9$ and $\rho^3 = 0.0$ are shown.
Although road $3$ is completely free, only few case can enter, as their way is blocked by cars waiting to enter road $2$. 
In this case $C^1=2 c^2$.
Thus the high density on road $2$ is causing a left going shock on the ingoing road.
A layer forms only on road $2$, since on the other two the macroscopic characteristics move away from the junction.

\section{Technical details}
\label{technical}

In this section we give the proof of Theorem \ref{theoremdiverge}.

\subsection{Coupling the layers (Diverging lanes with no driver preferences)}		
\label{layerproof}

In the first  step of the proof of Theorem  \ref{theoremdiverge}  the combination of the kinetic coupling conditions  with the layer equations has to be  considered.
Each layer can have either a stable solution (S) or an unstable solution (U).
Thus, for three edges we have eight possible combinations, which be denote by U/S-U/S-U/S.

\noindent{\bf Case 1, U-U-U.}	
We have $\rho^1 = \rho_+ (C^1), \rho^2= \rho_-(C^2) \rho^3= \rho_-(C^3)$ with $0 \le  C^1,C^2,C^3 < \sigma$.
We consider the four different cases for the kinetic coupling conditions from section \ref{kincondsummary}.

{\bf Case I:}\vspace{-1.5em}
\begin{align*}
\rho_+(C^1) &= \rho_-(C^2)= \rho_-(C^3)\\
C^1&= C^2+C^3.
\end{align*}	

The first two equalities give
$C^1=C^2=C^3 = \sigma$. This is not consistent with the range of $C^i$.

{\bf Case II:}\vspace{-1.5em}
\begin{align*}
\rho_+(C^1) &= C^1 \ge \rho_-(C^2) = \rho_-(C^3) \\
C^1&= C^2+C^3.
\end{align*}	
The first equality  has no solution, since $F(\rho) \le \rho$ and $F$  strictly concave gives 
$\sigma <\rho^\star =\rho_+(\sigma)$. 

{\bf Case III:}\vspace{-1.5em}
\begin{align*}
C^1 \le \rho_+(C^1) &= \rho_-(C^3)\le 
\rho_-(C^2)  \\
C^2&=0, C^1= C^2+C^3.
\end{align*}	
The first equality  gives $C^1=C^3 =\sigma$, which is not in the range of the $C^1,C^3$.

{\bf Case IV:}\vspace{-1.5em}
\begin{align*}
C^1 \le \rho_+(C^1) &= \rho_-(C^2)
\le  \rho_-(C^3)\\
C^3&=0, C^1= C^2+C^3.
\end{align*}	
The  first equality gives
$C^1=C^2=\sigma$,  which is not in the range of $C^1, C^2$. 

Altogether this  combination is not admissible.

\noindent{\bf Case 2, S-U-U} 

We have $\rho_0^1  \in [0,\rho_+(C^1)), \rho^2_0= \rho_-(C^2) , \rho^3_0= \rho_-(C^3)$ with
 $0 \le  C^1 \le  \sigma$, $0 \le  C^2,C^3 < \sigma  $.

{\bf Case I:}\vspace{-1.5em}
\begin{align*}
\rho_0^1 &= \rho_-(C^2)= \rho_-(C^3) \ge C^1\\
C^1&= C^2+C^3.
\end{align*}	
The equalities  give
$C^2= C^3= \frac{C^1}{2}$. Moreover, $\rho_0^1  \in [0,\rho_+(C^1))$ is obviously consistent with the above conditions. Finally, we have the condition $\rho_-(\frac{C^1}{2}) \ge C^1$.

{\bf Case II:}\vspace{-1.5em}
\begin{align*}
C^1 &= \rho_0^1  \ge 
\rho_-(C^2) = \rho_-(C^3) \\
C^1&= C^2+C^3.
\end{align*}	
The equalities give  $C^2=C^3=\frac{C^1}{2}$. 
$\rho_0^1  \in [0,\rho_+(C^1))$ is consistent with the conditions. We have finally $C^1 \ge \rho_-(\frac{C^1}{2})$.

{\bf Case III:} \vspace{-1.5em}
\begin{align*}
C^1 \le \rho_0^1 &= \rho_-(C^3)\le 
\rho_-(C^2)  \\
C^2&=0\;,\   C^1= C^2+C^3\ .
\end{align*}	
This leads to $0 \ge \rho_-(C^3) $ or $C^3=0$ and then $C^1=0$.

{\bf Case IV:}\vspace{-1.5em}
\begin{align*}
C^1 \le \rho_0^1  &= \rho_-(C^2) \le 
  \rho_-(C^3)\\
C^3&=0\;,\ C^1= C^2+C^3\ .
\end{align*}	
This leads to $0 \ge \rho_-(C^2) $ or $C^2=0$ and then $C^1=0$.

\noindent{\bf Case 3, U-S-U} 

Here it is $\rho^1 = \rho_+ (C^1)$, $\rho^2 \in (\rho_-(C^2),1)$, $\rho^3= \rho_-(C^3)$ with
 $0 \le  C^2 \le  \sigma$, $0 \le  C^1,C^3 < \sigma  $.

 {\bf Case I:}\vspace{-1.5em}
\begin{align*}
\rho_+(C^1) &= \rho^2_0 = \rho_-(C^3) \ge C^1\\
C^1&= C^2+C^3.
\end{align*}	
The equalities  give
$C^1= C^3=\sigma$ which is not in the range of $C^1,C^3$.

{\bf Case II:}\vspace{-1.5em}
\begin{align*}
C^1 &= \rho_+ (C^1) \ge 
\rho^2_0 = \rho_-(C^3) \\
C^1&= C^2+C^3.
\end{align*}	
The first  equation has no solution.

{\bf Case III:}\vspace{-1.5em}
\begin{align*}
C^1 \le \rho_+ (C^1)&= \rho_-(C^3) \le 
\rho^2_0 \\
C^2&=0\;,\ C^1= C^2+C^3\ .
\end{align*}	
This leads to $C^1 =C^3 = \sigma$ which is not in the range of $C^1,C^3$.

{\bf Case IV:}\vspace{-1.5em}
\begin{align*}
C^1 \le \rho_+ (C^1) &= \rho^2_0 \le 
 \rho_-(C^3)\\
C^3&=0\;,\  C^1= C^2+C^3\ .
\end{align*}	
$\rho^2_0 \in (\rho_-(C^2),1)$  leads to $0=\rho_-(C^3) \ge \rho_0^2  \ge \rho_-(C^2)$ or  $C^2=0$ and $C^1=0$.
This leads to a contradiction to $1=\rho_+(0) \le \rho_-(C^3) = 0$.

This case is not admissible.

\noindent{\bf Case 4, U-U-S} 

This case is symmetric to Case 3 and not admissible.

\noindent{\bf Case 5, U-S-S} 

We have $\rho^1 = \rho_+ (C^1), \rho^2 \in  (\rho_-(C^2),1), \rho^3 \in  (\rho_-(C^3),1)$
with  $0 \le  C^2,C^3 \le  \sigma$, $0 \le  C^1 < \sigma  $.

{\bf Case I:}\vspace{-1.5em}
\begin{align*}
C^1 \le \rho_+(C^1) &= \rho^2_0 
= \rho_0^3\\
C^1&= C^2+C^3\ .
\end{align*}	
The conditions are consistent  with the requirements for $\rho^2_0 $ and $\rho_0^3$. We need  $0 \le C^2+C^3 < \sigma$. One obtains $\rho_+(C^1) = \rho^1_0= \rho^2_0 
= \rho_0^3$.

{\bf Case II:}\vspace{-1.5em}
\begin{align*}
C^1 &= \rho_+ (C^1) \ge 
\rho^2_0= \rho_0^3\\
C^1&= C^2+C^3\ .
\end{align*}	
The first  equation has no solution.

{\bf Case III:}\vspace{-1.5em}
\begin{align*}
C^1 \le \rho_+ (C^1)&= \rho_0^3
\le \rho^2_0 \\
C^2&=0\;,\  C^1= C^2+C^3\ .
\end{align*}	
We have $0 \le C^1=C^3< \sigma $. 
We have $\rho^2_0 \ge \rho_+ (C^1)= \rho_0^1= \rho_0^3$.

{\bf Case IV:}\vspace{-1.5em}
\begin{align*}
C^1 \le \rho_+ (C^1) &= \rho^2_0 \le   \rho_0^3\\
C^3&=0\;,\  C^1= C^2+C^3\ .
\end{align*}	
The conditions are consistent  with the requirements for  $\rho_0^2$. We have $0 \le C^1=C^2< \sigma $. 
We have $\rho^3_0 \ge \rho_+ (C^1)= \rho_0^1= \rho_0^2$.

\noindent{\bf Case 6, S-U-S} 

We have $\rho_0^1 \in  [0,\rho_+ (C^1)), \rho_0^2= \rho_-(C^2) ,\rho_0^3 \in (\rho_-(C^3),1)$
with  $0 \le  C^1,C^3 \le  \sigma$, $0 \le  C^2 < \sigma  $.

{\bf Case I:}\vspace{-1.5em}
\begin{align*}
C^1 \le \rho^1 _0 &= \rho_-(C^2) 
= \rho_0^3\\
C^1&= C^2+C^3\ .
\end{align*}
The range of $\rho_0^3$ gives $\rho_-(C^2)  \ge \rho_-(C^3)$ or $C^2 \ge C^3$ or $C^2 \ge C^1-C^2$
or $C^1 \le 2 C^2$. In the same way,  $C^1 \ge 2 C^3$ and therefore $2 C^3 \le C^1 \le 2 C^2$.
Moreover, we need 
$\rho_-(C^2) \ge C^1$.

We have then  $\rho^1 _0=\rho^2_0 = \rho^3 _0=\rho_-(C^2)$.

{\bf Case II:}\vspace{-1.5em}
\begin{align*}
\rho^1 _0  &= C^1 \ge 
\rho_-(C^2)= \rho_0^3\\
C^1&= C^2+C^3\ .
\end{align*}	
The range of $\rho_0^3$ gives $\rho_-(C^2)  \ge \rho_-(C^3)$ or  $\rho_-(C^2)  \ge \rho_-(C^1-C^2)$  or  $C^1 \le 2 C^2 $ and similarly
$C^1\ge 2 C^3$. Together $2 C^3 \le C^1 \le 2 C^2 $. Moreover, we need $C^1 \ge 
 \rho_-(C^2)$.
 Then
 $\rho^1 _0=C^1 \ge \rho^2_0 = \rho^3 _0=\rho_-(C^2)$.

{\bf Case III:}\vspace{-1.5em}
\begin{align*}
C^1 \le \rho^1 _0 &= \rho_0^3 \le
\rho_-(C^2) \\
C^2&=0\;,\  C^1= C^2+C^3\ .
\end{align*}	
The ranges of $\rho^1_0$ and $\rho^3_0$ lead  to $\rho_+ (C^1) \ge \rho_-(C^3)$ and  $\rho_-(C^3)  \le \rho_-(C^2)=0$ which gives $C^3=0$
and $C^1=0$. Moreover, we have $\rho_0^1=0$  and then  $\rho_0^3=\rho_0^1=\rho_0^2=0$.

{\bf Case IV:}\vspace{-1.5em}
\begin{align*}
C^1 \le \rho^1 _0 &= \rho_-(C^2) 
\le  \rho_0^3\\
C^3&=0\;,\  C^1= C^2+C^3\ .
\end{align*}	
We have   $C^1=C^2< \sigma$. 
Moreover, $\rho_0^3 \ge \rho_0^2 =\rho_0^1=\rho_-(C^2)\ge C^1$.

\noindent{\bf Case 7, S-S-U}

The case is symmetric to case 6.

\noindent{\bf Case 8, S-S-S} 

We have $\rho_0^1  \in  [0,\rho_+ (C^1)), \rho^2 \in (\rho_-(C^2),1) ,\rho_0^3 \in ( \rho_-(C^3),1)$
with 

$0  \le C^1,C^2,C^3 \le \sigma$.

{\bf Case I:}\vspace{-1.5em}
\begin{align*}
C^1 \le \rho^1 _0 &= \rho^2_0
= \rho^3_0\\
C^1&= C^2+C^3\ .
\end{align*}	

There is no constraint  on the $C^i$ except $0 \le C^2 + C^3 \le \sigma$ .

We have $\max( C^1, \rho_-(C^2), \rho_-(C^3)) \le \rho_0^i  \le \rho_+(C^1)$

{\bf Case II:}\vspace{-1.5em}
\begin{align*}
\rho^1 _0  &= C^1 \ge 
\rho_0^2=\rho_0^3 \\
C^1&= C^2+C^3\ .
\end{align*}	
We have $\rho_-(C^2) , \rho_-(C^3)\le C^1$. Then $\rho^1 _0=C^1 \ge 
\rho^2 _0=\rho^3 _0 \ge \max(\rho_-(C^2), \rho_-(C^3))$.

{\bf Case III:}\vspace{-1.5em}
\begin{align*}
C^1 \le \rho^1 _0 &= \rho_0^3 \le 
\rho_0^2 \\
C^2&=0\;,\  C^1= C^2+C^3\ .
\end{align*}	
We have $C^1=C^3$ and 
$\rho_-(C^3)\le \rho_0^3=\rho_0^1\le \min(\rho_0^2, \rho_+(C^1))$.

{\bf Case IV:}\vspace{-1.5em}
\begin{align*}
C^1 \le \rho^1 _0 &= \rho_0^2
\le  \rho_0^3 \\
C^3&=0\;,\  C^1= C^2+C^3\ .
\end{align*}	
Then $C^1=C^2$ and
$\rho_+(C^1),\rho_0^3 \ge  \rho_0^2 =\rho_0^1 \ge \rho_-(C^2)$.

\subsubsection{Summary}
\label{layersum}
Altogether we obtain

\noindent{\bf Case1, U-U-U.}	
This  combination is not admissible.

\noindent{\bf Case 2, S-U-U} 
Only, if   $C^2=C^3=\frac{C^1}{2}$.

If $\rho_-(\frac{C^1}{2}) >C^1$, then, 
$\rho_0^2 = \rho_0^3=\rho_0^1 = \rho_-(\frac{C^1}{2}) $.

If  $\rho_-(\frac{C^1}{2}) <C^1$, then 
$\rho_0^1 =C^1$ and  $\rho_0^2 = \rho_0^3= \rho_-(\frac{C^1}{2}) $.

\noindent{\bf Case 3, U-S-U} This case is not admissible.

\noindent{\bf Case 4, U-U-S} This case is not admissible.

\noindent{\bf Case 5, U-S-S} 
If  $0 \le C^2+C^3 < \sigma$ and $C^2\neq0 \neq C^3 $ we have 
$\rho_+ (C^1)=\rho^1_0=\rho^2_0=\rho^3_0 $.

There are two special cases:
\begin{enumerate}
\item 
If $C^2=0$ and  $0 \le C^1=C^3 < \sigma$, then $\rho_0^2$ is not uniquely determined with \\
$\rho^2_0 \ge  \rho^1_0=\rho^3_0=\rho_+ (C^1) $.
\item
If $C^3=0$, $0 \le C^1,C^2 < \sigma$ then $\rho_0^3$ is not uniquely determined with \\
$\rho^3_0 \ge  \rho^1_0=\rho^2_0 = \rho_+ (C^1)$.
\end{enumerate}

\noindent{\bf Case 6, S-U-S} 
This case requires $2 C^3 \le C^1 \le 2 C^2 $, $C^2 < \sigma$.
We have   two subcases:
\begin{enumerate}
\item 
If  $\rho_-(C^2) \ge C^1$, $C^3 \neq0$, then 
  $\rho^1 _0=\rho^2_0 = \rho^3 _0=\rho_-(C^2)$.
In the special case  
 $C^3=0$ we have that   $\rho^3 _0$ is not uniquely determined with 
$\rho^3_0 \ge \rho_-(C^2) =\rho_0^2 =\rho_0^1\ge C^1$. 
\item  
 If $\rho_-(C^2) \le C^1$ one obtains
$\rho^1_0=C^1 $ 
and $\rho^2_0 = \rho^3 _0=\rho_-(C^2)$.
\end{enumerate}
\noindent{\bf Case 7, S-S-U}
The case is symmetric to case 6.

This case requires $2 C^2 \le C^1 \le 2 C^3 $, $C^3 < \sigma$.
We have   two subcases:
\begin{enumerate}
\item 
If  $\rho_-(C^3) \ge C^1$, $C^2 \neq0$, then 
  $\rho^1_0=\rho^2_0 = \rho^3 _0=\rho_-(C^3)$.
In the special case  
 $C^2=0$ we have that   $\rho^2 _0$ is not uniquely determined with 
$\rho^2_0 \ge \rho_-(C^3) =\rho_0^3 =\rho_0^1\ge C^1$. 
\item  
 If $\rho_-(C^3) \le C^1$ one obtains
$\rho^1_0=C^1 $ 
and $\rho^2_0 = \rho^3_0=\rho_-(C^3)$.
\end{enumerate}

\noindent{\bf Case 8, S-S-S} 
We need $0 \le C^2+C^3 \le \sigma$. $\rho_0^1,\rho_0^2, \rho_0^3$ are undetermined.

We have two subcases
\begin{enumerate} 
\item
If $C^1 \le \max(\rho_-(C^2), \rho_-(C^3)) $, $C^1\neq C^2$ and $C^1 \neq C^3$, then \\
 $\max( \rho_-(C^2), \rho_-(C^3)) \le \rho_0^3= \rho_0^2=\rho^1_0 \le  \rho_+(C^1)$.\\
We have two special cases:\\
For  $C^1=C^3$ one has  $\rho_-(C^3)\le \rho_0^3=\rho_0^1\le \min( \rho_+(C^1),\rho_0^2)$.\\
For $C^1=C^2$ one has $\rho_-(C^2) \le  \rho_0^2 =\rho_0^1 \le \min(\rho_+(C^1),\rho_0^3)$.
\item
For  $\max(\rho_-(C^2), \rho_-(C^3)) \le C^1$, one has \\
$ \max(\rho_-(C^2), \rho_-(C^3))\le \rho^2 _0=\rho^3 _0 \le \rho^1 _0=C^1$.
\end{enumerate}


\subsection{Matching layer solutions and half Riemann problems (Diverging lanes with no driver preferences)}
\label{proof}

Assuming  the  initial states $\rho_B^1,\rho_B^2,\rho_B^3$ to be given, we have to determine the fluxes $C^i$ and new states $\rho_K^i$ at the node. As mentioned, on the one hand $\rho_K^i$ are the asymptotic states of the respective layer problems and they and  the corresponding fluxes $C^i$ have to fulfill the conditions on the single kinetic layers, see section \ref{summary}, and on the coupled layers, see section \ref{layersum}.
On the other hand they are the left (road 1 and 2) or right hand (road 3) states of the half Riemann problems with $\rho_B^i$
as corresponding states fullfillling the conditions in \ref{Riemann}.
We  consider eight different configurations for the states $\rho_B^i$ corresponding to combinations of different half Riemann problems.
For  each of them all possible combinations with stable or unstable layer solutions have to be discussed.
Not admissible combinations are not listed. We use the  notation
$$
\rho_-^{D} (C) = \max (\rho_-(C), D), 0 \le C,D \le \sigma.
$$

\noindent{\bf Case 1, RP1-1-1} $\rho_B^1 \ge \rho^\star , \rho_B^2 \le \rho^\star , \rho_B^3 \le  \rho^\star $.
From section \ref{summary} we obtain 
\begin{align*}
\rho_K^1 &\in [\rho^\star,1] :
&  (U) &\text{ or } ((S) \text{ with } C^1=\sigma) \\
\rho_K^2 &\in [0, \rho^\star]: 
& (U) &\text{ or } ((S) \text{ with } C^2=\sigma)\\
\rho_K^3 &\in [0,\rho^\star]:
& (U) &\text{ or } ((S) \text{ with } C^3=\sigma)\ .
\end{align*}
The discussion in section 	\ref{layersum} gives the following 5 different cases:
\begin{enumerate}
	\item[{\bf SUU}] with $C^1 = \sigma$  and $C^2=C^3=\frac{\sigma}{2}$, $\rho_0^1=\rho_-^{C^1} (\frac{C^1}{2}), \rho_0^2=\rho_0^3=\rho_-(\frac{C^1}{2}) $.
	\item[{\bf USS}] with $C^2 = C^3 = \sigma$ 
	which contradicts the balance of fluxes.
	\item[{\bf SUS}] with $C^1 = C^3 = \sigma$ which contradicts  $2 C^3 \le  C^1$ and $C^1=C^2$.
	\item[{\bf SSU}] with $C^1 = C^2 = \sigma$ and a contradiction to $2 C^2 \le C^1$ and $C^1=C^3$. 
	\item[{\bf SSS}] with $C^1 = C^2 = C^3 = \sigma$, which gives a contradiction to the balance of fluxes.		
\end{enumerate}
This gives 
\begin{align*}
C^1 = \sigma\; \mbox{and} \;C^2=C^3=\frac{\sigma}{2},
\rho_0^1=\rho^\sigma_-(\frac{\sigma}{2}), \rho_0^2=\rho_0^3=\rho_-(\frac{\sigma}{2})\ .
\end{align*}
%

\noindent{\bf Case 2, RP1-1-2} $\rho_B^1 \ge \rho^\star , \rho_B^2 \le \rho^\star , \rho_B^3 \ge  \rho^\star $.
\begin{align*}
\rho_K^1 &\in [\rho^\star,1] :
&  (U) &\text{ or } ((S) \text{ with } C^1=\sigma) \\
\rho_K^2 &\in [0,\rho^\star]: 
& (U) &\text{ or } ((S) \text{ with } C^2=\sigma)\\		
\rho_K^3 &\in [0,\tau(\rho_B^3)] \cup \{\rho_B^3\}:
& ((U) \text{ with } C^3 \le F(\rho_B^3)) &\text{ or } ((S)
\text{ with } C^3=F(\rho_B^3))
\end{align*}	
\begin{enumerate}
	\item[{\bf SUU}] with $C^1 = \sigma$  and $C^2=C^3= \frac{\sigma}{2}$, if  $C^3 =\frac{\sigma}{2}\le F(\rho_B^3))$.
	Then  $\rho_0^1=\rho^{C^1}_-(\frac{C^1}{2}), \rho_0^2=\rho_0^3=\rho_-(\frac{C^1}{2}) $.
	\item[{\bf USS}] with $C^2 =  \sigma$. This contradicts $C^2+C^3< \sigma$.
	\item[{\bf SUS}] with $C^1 =  \sigma$,$C^3 = F(\rho_B^3)$ and $ F(\rho_B^3)\le \frac{\sigma}{2}$. Then  $\rho_0^1=\rho^{C^1}_-(C^2)$, $\rho_0^2=\rho_0^3=\rho_-(C^2) $.
	\item[{\bf SSU}] with $C^1 = C^2 = \sigma$ and  $0=C^3 $ and a contradiction to $C^1\le 2 C^3$ and $C^1=C^2 < \sigma$.
	\item[{\bf SSS}] with $C^1 = C^2 =  \sigma, C^3=F(\rho_B^3)=0$. This is  only possible, if $\rho_B^3=1$. 	Then $\rho_0^1=\rho_0^2=\rho_0^3=\rho^\star $.	
\end{enumerate}

This gives for 
\begin{align*}
F(\rho_B^3) \le \frac{\sigma}{2}&: C^1 =  \sigma, C^3 = F(\rho_B^3),\\& \rho_0^1=\rho^\sigma_-(\sigma - F(\rho_B^3)), \rho_0^2=\rho_0^3=\rho_-(\sigma - F(\rho_B^3)) \\
\frac{\sigma}{2}\le F(\rho_B^3)&: C^1 = \sigma, C^2=C^3= \frac{\sigma}{2}, \rho_0^1=\rho^\sigma_-(\frac{\sigma}{2}), \rho_0^2=\rho_0^3=\rho_-(\frac{\sigma}{2}).
\end{align*}

\noindent{\bf Case 3, RP1-2-1} $\rho_B^1 \ge \rho^\star , \rho_B^2 \ge \rho^\star , \rho_B^3 \le  \rho^\star $. 
Symmetric to Case 2.
 \begin{align*}
 \frac{\sigma}{2}\ge F(\rho_B^2))&: C^1 =  \sigma, C^2 = F(\rho_B^2),\\& \rho_0^1=\rho^\sigma_-(\sigma - F(\rho_B^2)), \rho_0^2=\rho_0^3=\rho_-(\sigma - F(\rho_B^2)) \\
 \frac{\sigma}{2}\le F(\rho_B^2))&: C^1 = \sigma, C^2=C^3= \frac{\sigma}{2}, \rho_0^1=\rho^\sigma_-(\frac{\sigma}{2}), \rho_0^2=\rho_0^3=\rho_-(\frac{\sigma}{2}).
 \end{align*}

\noindent{\bf Case 4, RP2-1-1} $\rho_B^1 \le \rho^\star , \rho_B^2 \le \rho^\star , \rho_B^3 \le  \rho^\star $. 
\begin{align*}
\rho_K^1 &\in  [1-\rho_B^1,1] \cup \{\rho_B^1\}:   
&  ((U)\text{ with } C^1 \le F(\rho_B^1)) &\text{ or } ((S) \text{ with } C^1= F(\rho_B^1))\\
\rho_K^2 &\in  [0,\rho^\star] : 
& (U)  &\text{ or } ((S)
\text{ with } C^2=\sigma)	\\	
\rho_K^3 &\in [0,1/2]:
& (U) &\text{ or } ((S) \text{ with } C^3=\sigma)	
\end{align*}	
\begin{enumerate}
	\item[{\bf SUU}] with $C^1 = F(\rho_B^1)$  and  $C^2= C^3 =\frac{ F(\rho_B^1)}{2}$. Moreover,  $\rho_0^1=\rho^{F(\rho_B^1)}_-(\frac{F(\rho_B^1)}{2})$, $\rho_0^2=\rho_0^3=\rho_-(\frac{F(\rho_B^1)}{2})$.
	\item[{\bf USS}] with $C^2 =  C^3=\sigma$ and a contradiction to the balance of fluxes.
	\item[{\bf SUS}] with $C^1 = F(\rho_B^1) , C^3 =\sigma$. $2 C^3 \le  C^1$ gives  a contradiction.
	\item[{\bf SSU}] with $C^1  = F(\rho_B^1)$ and $C^2 = \sigma$. $2 C^2 \le  C^1$ gives  a contradiction.
	\item[{\bf SSS}] with $C^1  = F(\rho_B^1)$, $C^2 =  C^3 = \sigma$ and a contradiction to the balance of fluxes.
\end{enumerate}
This gives alltogether 
\begin{align*}
  C^2= C^3 =\frac{C^1}{2} =\frac{ F(\rho_B^1)}{2}\ ,\quad
\rho_0^1=\rho^{F(\rho_B^1)}_-(\frac{F(\rho_B^1)}{2})\ ,\quad  \rho_0^2=\rho_0^3=\rho_-(\frac{F(\rho_B^1)}{2})\ .
\end{align*}

\noindent{\bf Case 5, RP1-2-2} $\rho_B^1 \ge \rho^\star , \rho_B^2 \ge \rho^\star , \rho_B^3 \ge  \rho^\star $. 
\begin{align*}
\rho_K^1 &\in  [\rho^\star,1] : 
& (U)  &\text{ or } ((S)
\text{ with } C^1=\sigma)	\\
\rho_K^2 &\in [0,1-\rho_B^2] \cup \{\rho_B^2\}:  
&  ((U)\text{ with } C^2 \le F(\rho_B^2)) &\text{ or } ((S) \text{ with } C^2= F(\rho_B^2))\\
\rho_K^3 &\in [0,1-\rho_B^3] \cup \{\rho_B^3\}: 
&((U)\text{ with } C^3 \le F(\rho_B^3)) &\text{ or } ((S) \text{ with } C^3= F(\rho_B^3))
\end{align*}	
\begin{enumerate}
	\item[{\bf SUU}] with $C^1 = \sigma$,  $C^2 \le  F(\rho_B^2)$  and $C^3 \le  F(\rho_B^3)$.
	With $C^2=C^3 =\frac{C^1}{2}$ one obtains $F(\rho_B^2)+  F(\rho_B^3) \ge \sigma$ and $F(\rho_B^2) \ge \frac{\sigma}{2}$, $F(\rho_B^3) \ge \frac{\sigma}{2}$. We have $\rho_0^1= \rho_-^\sigma(\frac{\sigma}{2}), \rho_0^2=\rho_0^3=\rho_-(\frac{\sigma}{2})$.
	\item[{\bf USS}] with $C^2 = F(\rho_B^2) , C^3 = F(\rho_B^3) $. With $C^2 + C^3\le \sigma$ or $F(\rho_B^2)+F(\rho_B^3) \le \sigma$,
	$\rho_0^1=\rho_0^2=\rho_0^3=\rho_+(F(\rho_B^2) + F(\rho_B^3) )$
	\item[{\bf SUS}] with $C^1 = \sigma , C^2 \le F(\rho_B^2), C^3 = F(\rho_B^3)$. Then   $2C^3 \le  C^1$
	gives $ F(\rho_B^3) \le \frac{\sigma}{2}$. Moreover, $F(\rho_B^2)+  F(\rho_B^3) \ge \sigma$.
	We have $\rho_0^1= \rho^\sigma_-(\sigma - F(\rho_B^3)), \rho_0^2=\rho_0^3=\rho_-(\sigma - F(\rho_B^3))$
	\item[{\bf SSU}] with $C^1  = \sigma$ and $C^2 = F(\rho_B^2),C^3 \le  F(\rho_B^3) $. $2 C^2 \le C^1$ gives
	$F(\rho_B^2) \le \frac{\sigma}{2}$. Moreover, $F(\rho_B^2)+  F(\rho_B^3) \ge \sigma$.
	$\rho_0^1= \rho^\sigma_-(\sigma - F(\rho_B^2)), \rho_0^2=\rho_0^3=\rho_-(\sigma - F(\rho_B^2))$
	\item[{\bf SSS}] with $C^1 = \sigma$ and $C^2 =F(\rho_B^2) , C^3 =F(\rho_B^3) $. This is only possible, if $F(\rho_B^2) + F(\rho_B^3) = \sigma$. 
	We have $\max(\rho_-(C^2),\rho_-(C^3) )\le \rho_0^2=\rho_0^3 =  \rho_0^1 \le  \rho_+(\sigma) = \rho^\star$
	for
	$\max(\rho_-(C^2), \rho_-(C^3)) \ge \sigma$ or	
	$ \max(\rho_-(C^2), \rho_-(C^3))\le \rho^2 _0=\rho^3 _0 \le \rho^1 _0=\sigma$
	for $\max(\rho_-(C^2), \rho_-(C^3)) \le \sigma$.
\end{enumerate}
This gives altogether for 
\begin{align*}
F(\rho_B^2)+  F(\rho_B^3)& \le \sigma:
C^2 = F(\rho_B^2) , C^3 = F(\rho_B^3)\ ,  \\
&\rho_0^1=\rho_0^2=\rho_0^3=\rho_+(F(\rho_B^2)+F(\rho_B^3))\\
F(\rho_B^2) \ge \frac{\sigma}{2},F(\rho_B^3) &\ge \frac{\sigma}{2}:	
C^1 = \sigma, C^2= \frac{\sigma}{2}, C^3= \frac{\sigma}{2}\ ,\\
&\rho_0^1= \rho^\sigma_-(\frac{\sigma}{2}), \rho_0^2=\rho_0^3=\rho_-(\frac{\sigma}{2})\ ,\\
F(\rho_B^2)+  F(\rho_B^3) \ge \sigma, F(\rho_B^3) &\le \frac{\sigma}{2}:C^1 = \sigma , C^2  = \sigma- F(\rho_B^3), C^3 = F(\rho_B^3)\ , \\
&\rho_0^1= \rho^\sigma_-(\sigma - F(\rho_B^3)), \rho_0^2=\rho_0^3=\rho_-(\sigma - F(\rho_B^3))\ ,\\ 
F(\rho_B^2)+  F(\rho_B^3) \ge \sigma, F(\rho_B^2) &\le \frac{\sigma}{2}:C^1  = \sigma,C^2 = F(\rho_B^2),C^3 =\sigma- F(\rho_B^2)\ ,\\
&\rho_0^1= \rho^\sigma_-(\sigma - F(\rho_B^2)), \rho_0^2=\rho_0^3=\rho_-(\sigma - F(\rho_B^2))\ .
\end{align*}
\begin{remark}
We note that at the interfaces between the different conditions we obtain values  for the $\rho_0^i $ which correspond to the range of values   for the $\rho_0^i$-values in case (SSS).
\end{remark}

\noindent{\bf Case 6, RP2-1-2} $\rho_B^1 \le \rho^\star , \rho_B^2 \le \rho^\star , \rho_B^3 \ge  \rho^\star $. 
\begin{align*}
\rho_K^1 &\in [1-\rho_B^1,1] \cup \{\rho_B^1\}:  
&  ((U)\text{ with } C^1 \le F(\rho_B^1)) &\text{ or } ((S) \text{ with } C^1= F(\rho_B^1))\\
\rho_K^2 &\in  [0,\rho^\star] : 
& (U)  &\text{ or } ((S)
\text{ with } C^2=\sigma)	\\
\rho_K^3 &\in [0,1-\rho_B^3] \cup \{\rho_B^3\}: 
&((U)\text{ with } C^3 \le F(\rho_B^3)) &\text{ or } ((S) \text{ with } C^3= F(\rho_B^3))
\end{align*}	
\begin{enumerate}
	\item[{\bf SUU}] with $C^1 =  F(\rho_B^1)$, $C^3  \le  F(\rho_B^3)$. $C^2=C^3= \frac{C^1}{2} = \frac{F(\rho_B^1)}{2}$ leads to
	$ F(\rho_B^1)  \le 2  F(\rho_B^3)$, $\rho_0^1=\rho^{F(\rho_B^1)}_-(\frac{F(\rho_B^1)}{2}), \rho_0^2=\rho_0^3=\rho_-(\frac{F(\rho_B^1)}{2}).$
	\item[{\bf USS}] with $C^1 = F(\rho_B^1) , C^2 =\sigma, C^3 = F(\rho_B^3) =0$ contradicts $C^2 < \sigma$.
	\item[{\bf SUS}] with $C^1 =  F(\rho_B^1), C^3 = F(\rho_B^3)$.  $2 C^3 \le C^1$ gives $ 2 F(\rho_B^3)\le   F(\rho_B^1)$.
	$\rho_0^1=\rho^{F(\rho_B^1)}_-(F(\rho_B^1)-F(\rho_B^3)), \rho_0^2=\rho_0^3=\rho_-(F(\rho_B^1)-F(\rho_B^3)).$ 
	In the special case $C^3=0$ we have $\rho_0^3\ge \rho_-(F(\rho_B^1)-F(\rho_B^3))=\rho_0^1=\rho_0^2$.
	\item[{\bf SSU}] with $C^1  = F(\rho_B^1)$ and $C^2 = \sigma,C^3 \le  F(\rho_B^3) $. This requires $C^3=0$ which contradicts $C^2 \le C^3$.
	\item[{\bf SSS}] with $C^1 = F(\rho_B^1) $ and $C^2 =\sigma, C^3 =F(\rho_B^3) $. Only possible if $C^1 =\sigma, C^3 = 0$.
	 $\rho_0^1=\rho_0^2=\rho^\star\le  \rho_0^3. $ 
\end{enumerate}
This gives altogether for $F(\rho_B^3)>0$ and 
\begin{align*}
 F(\rho_B^1) & \le 2  F(\rho_B^3):
C^1 =  F(\rho_B^1),C^2=C^3= \frac{C^1}{2} = \frac{F(\rho_B^1)}{2}\ ,\\
&\rho_0^1= \rho^{F(\rho_B^1)}_-(\frac{F(\rho_B^1)}{2}), \rho_0^2=\rho_0^3=\rho_-(\frac{F(\rho_B^1)}{2})\ ,\\
2 F(\rho_B^3)&\le   F(\rho_B^1):	
C^1 =  F(\rho_B^1), C^3 = F(\rho_B^3), C^2=F(\rho_B^1)- F(\rho_B^3)\ ,\\
&\rho_0^1=\rho^{F(\rho_B^1)}_-(F(\rho_B^1)- F(\rho_B^3)), \rho_0^2=\rho_0^3=\rho_-(F(\rho_B^1)- F(\rho_B^3))\ .
\end{align*}
In the special case $F(\rho_B^3)=0$ the value of  $\rho_0^3$ is not uniquely determined.

\noindent{\bf Case 7, RP2-2-1} $\rho_B^1 \le \rho^\star , \rho_B^2 \ge \rho^\star , \rho_B^3 \le  \rho^\star $. 
Symmetric to Case 6. For $F(\rho_B^2)>0$ and 
\begin{align*}
F(\rho_B^1)  &\le 2  F(\rho_B^2):
C^1 =  F(\rho_B^1),C^2=C^3= \frac{C^1}{2} = \frac{F(\rho_B^1)}{2}.\\
&\rho_0^1= \rho^{F(\rho_B^1)}_-(\frac{F(\rho_B^1)}{2}), \rho_0^2=\rho_0^3=\rho_-(\frac{F(\rho_B^1)}{2}).\\
2 F(\rho_B^2)&\le   F(\rho_B^1):	
C^1 =  F(\rho_B^1), C^2 = F(\rho_B^2), C^3=F(\rho_B^1)- F(\rho_B^2).\\
&\rho_0^1=\rho^{F(\rho_B^1)}_-(F(\rho_B^1)- F(\rho_B^2)), \rho_0^2=\rho_0^3=\rho_-(F(\rho_B^1)- F(\rho_B^2)).
\end{align*}
In the special case $F(\rho_B^2)=0$ the value of  $\rho_0^2$ is not uniquely determined.

\noindent{\bf Case 8, RP2-2-2} $\rho_B^1 \le \rho^\star , \rho_B^2 \ge \rho^\star , \rho_B^3 \ge  \rho^\star $.
\begin{align*}
\rho_K^1 &\in [1-\rho_B^1,1] \cup \{\rho_B^1\}:  
&  ((U)\text{ with } C^1 \le F(\rho_B^1)) &\text{ or } ((S) \text{ with } C^1= F(\rho_B^1))\\
\rho_K^2 &\in [0,1-\rho_B^2] \cup \{\rho_B^2\}:  
&  ((U)\text{ with } C^2 \le F(\rho_B^2)) &\text{ or } ((S) \text{ with } C^2= F(\rho_B^2))\\
\rho_K^3 &\in [0,1-\rho_B^3] \cup \{\rho_B^3\}: 
&((U)\text{ with } C^3 \le F(\rho_B^3)) &\text{ or } ((S) \text{ with } C^3= F(\rho_B^3))
\end{align*}
\begin{enumerate}
	\item[{\bf SUU}] with $C^1 =   F(\rho_B^1)$, $C^2 \le   F(\rho_B^2)$. $C^3 \le   F(\rho_B^3)$.
	$C^2=C^3 =\frac{C^1}{2}$  gives $C^2=C^3 =\frac{F(\rho_B^1)}{2}$. We need
	$ F(\rho_B^1)\le   2 F(\rho_B^2),F(\rho_B^1)\le   2F(\rho_B^3)$.
	$\rho_0^1= \rho^{F(\rho_B^1)}_-(\frac{F(\rho_B^1)}{2}), \rho_0^2=\rho_0^3=\rho_-(\frac{F(\rho_B^1)}{2}).$
	\item[{\bf USS}] with $C^1 \le  F(\rho_B^1) , C^2 =F(\rho_B^2), C^3 = F(\rho_B^3) $. With $C^2+C^3\le \sigma$ we have 
	$ F(\rho_B^2) +   F(\rho_B^3)\le \sigma$.
	$\rho_0^1=\rho_+(F(\rho_B^1))=\rho_0^2=\rho_0^3.$
	\item[{\bf SUS}] with $C^1 =  F(\rho_B^1), C^2\le  F(\rho_B^2),C^3 = F(\rho_B^3)$.  $2 C^3 \le  C^1$ gives $ 2 F(\rho_B^3)\le  F(\rho_B^1)$ and $F(\rho_B^3)- F(\rho_B^1) \le F(\rho_B^2)$ or  $F(\rho_B^1)+ F(\rho_B^2) \ge F(\rho_B^3)$. 
	$\rho_0^1=\rho^{F(\rho_B^1)}_-(F(\rho_B^2)),\rho_0^2=\rho_0^3=\rho_-(F(\rho_B^2))$
	\item[{\bf SSU}] with $C^1  = F(\rho_B^1)$ and $C^2 =  F(\rho_B^2) ,C^3 \le  F(\rho_B^3) $. $2C^2 \le C^1$ gives
	$2F(\rho_B^2) \le F(\rho_B^1)$. $\rho_0^1=\rho^{F(\rho_B^1)}_-(F(\rho_B^3)),\rho_0^2=\rho_0^3=\rho_-(F(\rho_B^3))$
	\item[{\bf SSS}] with $C^1 = F(\rho_B^1) $ and $C^2 =F(\rho_B^2), C^3 =F(\rho_B^3) $. Only possible if $F(\rho_B^1)=F(\rho_B^2) + F(\rho_B^3)$. 	We have $F(\rho_B^2) + F(\rho_B^3) \le \max(\rho_-(F(\rho_B^2)),\rho_-(F(\rho_B^3) )\le \rho_0^2=\rho_0^3 =  \rho_0^1 \le  \rho_+(F(\rho_B^2) + F(\rho_B^3))$
		or	
		$ \max(\rho_-(F(\rho_B^2)), \rho_-(F(\rho_B^3))\le \rho^2 _0=\rho^3 _0 \le \rho^1 _0=F(\rho_B^2) + F(\rho_B^3)$.
\end{enumerate}
Additionally, there are again special cases with $F(\rho_B^2)=0$ or $F(\rho_B^3)=0$.

This gives altogether for $F(\rho_B^2)>0$ and  $F(\rho_B^3)>0$
\begin{align*}
 F(\rho_B^2) +   F(\rho_B^3)&\le F(\rho_B^1) 
 :
C^1 = F(\rho_B^2)+ F(\rho_B^3) , C^2 =F(\rho_B^2)\ , \\&C^3 = F(\rho_B^3),
\rho_0^1=\rho_+(F(\rho_B^1))=\rho_0^2=\rho_0^3\ ,\\
F(\rho_B^1)\le 2  &\min(F(\rho_B^2),    F(\rho_B^3) )
:
C^2= C^3= \frac{F(\rho_B^1)}{2}, C^1 = F(\rho_B^1)\ ,\\
&\rho_0^1= \rho^{F(\rho_B^1)}_-(\frac{F(\rho_B^1)}{2}), \rho_0^2=\rho_0^3=\rho_-(\frac{F(\rho_B^1)}{2})\ ,\\
F(\rho_B^2)+ F(\rho_B^3) &\ge F(\rho_B^1):\\
2 F(\rho_B^3)\le  & F(\rho_B^1) 
:
C^1 =  F(\rho_B^1), C^2= F(\rho_B^1) -F(\rho_B^3),C^3 = F(\rho_B^3)\ ,\\
&\rho_0^1=\rho^{F(\rho_B^1)}_-(F(\rho_B^2)),\rho_0^2=\rho_0^3=\rho_-(F(\rho_B^2))\ ,\\
2F(\rho_B^2)\le &  F(\rho_B^1) 
:
C^1  = F(\rho_B^1), C^2 =  F(\rho_B^2) ,C^3  =F(\rho_B^1)- F(\rho_B^2),\\ &\rho_0^1=\rho^{F(\rho_B^1)}_-(F(\rho_B^3)),\rho_0^2=\rho_0^3=\rho_-(F(\rho_B^3))\ .
\end{align*}
\begin{remark}
We note that at the interfaces between the different conditions we obtain values  for the $\rho_0^i $ which correspond to the range of values   for the $\rho_0^i$-values in case (SSS).
\end{remark}


\section{Conclusions}

	We have introduced  coupling conditions  for junctions with diverging lanes for a  kinetic  two velocity traffic model, which is used as a relaxation model
	for scalar traffic flow equations.
	From these coupling conditions we have derived, via  asymptotic analysis of the spatial layers at the nodes and a detailed investigation of the associated Riemann problems,  coupling conditions for classical scalar macroscopic traffic models.	
	The  derivation shows that a classical  condition  for a nonlinear scalar conservation law
	can be interpreted on the kinetic level as a combination of the balance of fluxes and a suitably modified equality of densites on all roads.
	This research is a continuation of the work  in  \cite{BK20}, where the case of merging lanes has been treated.

\end{document}